\documentclass[preprint,12pt]{elsarticle}
\usepackage{comment}
\usepackage{algorithm,algorithmic}
\usepackage{caption}
\usepackage{subcaption}
\usepackage[short]{optidef}
\usepackage{float}
\usepackage{adjustbox}
\usepackage{multicol}
\usepackage{multirow}
\usepackage{amssymb}
\usepackage{hyperref}
\hypersetup{
    colorlinks,
    citecolor=blue,
    linkcolor=blue,
}
\usepackage{amsmath}

\usepackage{geometry}
\newgeometry{left=2cm,right=2cm,top=2cm,bottom=2cm}

\begin{document}

\begin{frontmatter}

\title{Grey wolf optimizer and whale optimization algorithm for stochastic inventory management of reusable products in a two-level supply chain}

\author[inst1]{Amir Hossein Sadeghi}
\author[inst2]{Erfan Amani Bani}
\author[inst2]{Ali Fallahi\corref{mycorrespondingauthor}}\cortext[mycorrespondingauthor]{Corresponding author}\ead{Ali.fallahi@ie.sharif.edu}

\affiliation[inst1]{organization={Department of Industrial and Systems Engineering},
            addressline={North Carolina State University}, 
            city={Raleigh},
            state={NC},
            country={USA}}
\affiliation[inst2]{organization={Department of Industrial Engineering},
            addressline={Sharif University of Technology}, 
            city={Tehran},
            country={Iran}}

\begin{abstract}
Product reuse and recovery is an efficient tool that helps companies to simultaneously address economic and environmental dimensions of sustainability. This paper presents a novel problem for stock management of reusable products in a single-vendor, multi-product, multi-retailer network. Several constraints, such as the maximum budget, storage capacity, number of orders, etc., are considered in their stochastic form to provide a more realistic framework. The presented problem is formulated as a constrained nonlinear mathematical model. The chance-constrained programming method is suggested to deal with the constraints' uncertainty. Regarding the nonlinearity of the model, grey wolf optimizer (GWO) and whale optimization algorithm (WOA) as two novel metaheuristics are presented as solution approaches, and the sequential quadratic programming (SQP) exact algorithm validates their performance. The parameters of algorithms are calibrated using the Taguchi method for the design of experiments. Extensive analysis is established by solving several numerical results in different sizes and utilizing several comparison measures. Also, the results are compared statistically using proper parametric and non-parametric tests. The analysis of the results shows a significant difference between the algorithms, and GWO has a better performance for solving the presented problem. In addition, both algorithms perform well in searching the solution space, where the GWO and WOA differences with the optimal solution of the SQP algorithm are negligible.
\end{abstract}



\end{frontmatter}
\section{Introduction} \label{sec:intro}
Companies' perspectives have changed from classic business principles to contemporary ideas like the supply chain in today's competitive environment. More specifically, the managers try to integrate the activities and processes of their supply chain to improve the overall performance of their companies \cite{chang2016supply,taghiyeh2020multi}. This change is obvious in various sectors e.g., manufacturing, retail, healthcare, etc \cite{fakhrabad2023impact,nikoubin2023relax,sahebi2023evaluating,hasan2023multi}. Integration of the supply chain necessitates collaboration across its entities, as well as coordination of information and material flows. Such  coordination can play a significant role in cost-cutting and improving value for customers \cite{braglia2003modelling}. Inventory planning and control is one of the highly important problems in supply chains that need special coordination between supply chain entities. Unpleasant or variable inventory in a supply chain causes the Bullwhip effect and double marginalization, which ultimately degrade the performance of the supply chain and may even contribute to the demise of businesses \cite{disney2003effect}. Therefore, operations research experts investigated the problem from several years ago.  The history of inventory management dates back to a century ago when Harris \cite{harris1990many} developed the classical economic order quantity (EOQ) model. Afterward, this basic inventory model was extended for a wide variety of products by considering several realistic assumptions, with the details that will be discussed in the literature review section.

Despite the appearance of broad literature, the inventory management of reusable products, as an important type of items, remained unaddressed for several years. Reusable products are products that have the potential for reuse and recovery after consumption. Surgical tools, car parts, copiers, and other products are a few examples of these items. In addition to economic benefits, the recovery and reuse of reusable products help companies address the environmental dimension of sustainability appropriately. Recently, the single-product and multi-product EOQ models for inventory management of reusable items were by Mokhtari \cite{mokhtari2018joint}, and Fallahi et al. \cite{fallahi2022constrained}, respectively.  In these works, one of the main assumptions is that the system is single-level, and both models ignore vendor costs in the decision-making process. In other words, there is no coordination between the upstream and downstream of the supply chain, and it is assumed that the retailers are solely responsible for decision-making. On the other hand, both works modeled and optimized their problem by considering deterministic parameters. However, several uncertainty sources may impact the parameters of inventory systems, and the deterministic formulation of the problem may negatively impact the performance \cite{pasandideh2015optimization,khalilpourazari2019robust}.

In this research, we try to bridge the research gap and develop a single-vendor, multi-product, multi-retailer problem for inventory management of reusable items in a two-level supply chain. Also, we try to establish a more realistic problem by considering the system's operational constraints under uncertainty in the availability of resources. A chance-constrained programming approach is suggested to handle this uncertainty. Due to the nonlinearity and the dimension of the developed constrained nonlinear model, the grey wolf optimizer (GWO) and whale optimization algorithm (WOA) are designed and implemented to solve the problem. The exact sequential quadratic programming approach is proposed (SQP) for the validation of the metaheuristics. The primary contributions of this study are as follows:
\begin{itemize}
    \item Developing a problem for stock management of reusable products in a single-vendor multi-product multi-retailer supply chain under uncertainty in operations constraints.
    \item Presenting a chance-constrained programming approach to address the stochastic constraints.
    \item Designing GWO metaheuristic algorithm as the solution approach.
    \item Designing WOA metaheuristic algorithm as the solution approach.
    \item Presenting the SQP exact algorithm for performance validation of metaheuristics.
\end{itemize}

The rest of this article is organized as follows. Section~\ref{sec:lit} provides a review of the relevant works in the literature. Section~\ref{sec:problem} provides the problem definition and the model formulation. Section~\ref{sec:sol} describes the presented solution methodologies to solve the developed model. Section~\ref{sec:result} provides the analysis of results to show the presented model's applicability and compare the algorithms' performance. Finally, Section~\ref{sec:con} provides a conclusion of the research and suggestions for future works.

\section{Related work} \label{sec:lit}
In this section, we will review the relevant papers in the inventory management literature. As mentioned before, Harris \cite{harris1990many} introduced the first optimization model, EOQ, to determine the optimal ordering policy in an inventory system. This model was developed for a single-level system and included several other simple assumptions, such as infinite replenishment rate, availability of resources such as budget and space, etc. Therefore, the authors tried to extend EOQ to bring the model into a real-world setting as much as possible. A while after EOQ, Taft \cite{taft1918most} developed the economic production quantity (EPQ) and relaxed the infinite replenishment rate assumption to determine the optimal production quantity for inventory management of manufacturing companies. After that, EOQ and EPQ models were extended by considering several assumptions and features. For example, previous authors developed EOQ and EPQ models for special types of products such as deteriorating items \cite{covert1973eoq}, substitutable items \cite{mokhtari2022economic,drezner1995eoq}, growing items \cite{rezaei2014economic}, etc. Also, researchers considered several realistic assumptions, such as preventive maintenance \cite{mokhtari2020extended}, investment \cite{porteus1985investing}, trade credits \cite{taleizadeh2013eoq,chung2003optimal}, sustainability concerns \cite{chen2013carbon,taleizadeh2018sustainable}, discount, process reliability \cite{tripathy2003eoq}, pricing decisions \cite{cheng1990eoq}, marketing policies \cite{subramanyam1981eoq}, presence of imperfect items \cite{salameh2000economic}, inflation and time-value of money \cite{chandra1985effects}, transportation policies \cite{fallahi2021sustainable}, inspection errors \cite{khan2011economic}, etc. 

However, all these models were single-level and did not consider the vendors in the decision-making framework. Managers understood the importance of coordination in the supply chain and tried to provide integrated models for decision-making the so that the entire supply chain performance is optimized. For example, Pasandideh et al. \cite{pasandideh2011genetic} presented the formulation of EOQ for a single-vendor single-retailer supply chain. The backorder shortage was allowed, and several operational constraints were incorporated. The complexity of the problem prompted authors to use the genetic algorithm (GA) as the solution approach.  In a similar work, Pasandideh et al. \cite{pasandideh2014optimization} extended the EPQ with shortage for a constrained single-vendor single-retailer system and solved the problem via GA. Also, Pasandideh et al. \cite{pasandideh2011parameter} also worked on a stochastic inventory model for a single-vendor multi-retailer system. The authors utilized an expected value approach to handle the uncertainty. Taleizadeh et al. \cite{taleizadeh2011multiple} investigated the inventory management of a multi-vendor multi-retailer supply chain with order size-dependent lead time and partial back. The demand was assumed to be uniformly distributed, and the harmony search (HS) metaheuristic was implemented to find the solutions. Chen et al. \cite{chen2014retailer} considered a delay in payments as a trade credit option in a two-level supply chain and derived the optimal ordering policy of the system under this option.

In some other works, Cárdenas-Barrón \& Sana \cite{cardenas2015multi} provided an EOQ model for a single-vendor single-retailer supply chain in which the demand depended on the promotional effort. Also, the retailer's delay in payment was possible in the considered supply chain. Khan et al. \cite{khan2017learning} focused on developing an EOQ model for a two-level supply chain with the production of defective products. They assumed that the inspection process is subjected to error and that the production time depends on learning. Tiwari et al. \cite{tiwari2018joint} studied a supply chain model of deteriorating items with pricing and inventory decisions. It was assumed that there is a partial trade credit contract for both levels of the network. Karimian et al. \cite{karimian2020geometric} employed a geometric programming method for uncertainty inventory management in a single-vendor multi-retailer supply chain. The problem's applicability was shown by solving the problem for a case study in the Iranian furniture supply chain. Pourmohammad-Zia et al. \cite{pourmohammad2021food} presented a new model for coordinating vendor and retailer inventory in a growing products supply chain using (vendor-managed inventory) VMI and a cost-sharing contract. Pourmohammad-Zia et al. \cite{pourmohammad2021dynamic} also aimed to determine pricing and replenishment policies of the growing items supply chain in another work. Recently, Asadkhani et al. \cite{asadkhani2022sustainable} presented a sustainable supply chain under some emission reduction regulations and the VMI-consignment stock (VMI-CS) contract. They considered repair, salvage, and disposal as potential options to deal with imperfect items.

The focus of this paper is on the inventory management of reusable items. For the first time, Mokhtari \cite{mokhtari2018joint} designed a new single-product EOQ problem for stock management of reusable products. The author assumed that the resources are infinitely available and solved the unconstrained model by determining the optimal order and recovery quantity of reusable products through an analytical derivative-based method. Recently, Fallahi et al. \cite{fallahi2022constrained} pointed out that this model is not practical for systems that deal with multiple products and limitations of resources. Consequently, they presented a multi-product extension of the previous work, and considered the limitations on the maximum budget and the storage capacity for usable and recoverable items. They solved the model using differential evolution (DE) and particle swarm optimization (PSO) metaheuristics, and also developed two new versions of these algorithms using an intelligent machine learning algorithm. Table~\ref{tab:lit} compares the novelties of our research against the features of past papers in the literature

\begin{table}[H]
\caption{The features of the relevant past works}
\label{tab:lit}
\centering
\begin{adjustbox}{width=\textwidth}
\begin{tabular}{llllllllllll l l}
\hline
\multirow{2}{*}{Article} & \multirow{2}{*}{Year} & \multirow{2}{*}{Product Type} & \multirow{2}{*}{Model Type} & \multicolumn{2}{l}{Vendor Number} & \multicolumn{2}{l}{Retailer Number} & \multicolumn{2}{l}{Product Number} & \multirow{2}{*}{\begin{tabular}[c]{@{}l@{}}Operational\\ Constraints\end{tabular}} & \multirow{2}{*}{Uncertainty} & \multicolumn{2}{c}{Solution Approach} \\ \cline{5-10} \cline{13-14} 
 &  &  &  & Single & Multiple & Single & Multiple & Single & Multiple &  &  & Exact & (Meta) heuristic* \\ \hline
Drezner et al. \cite{drezner1995eoq} & 1995 & Substitutable & EOQ &  &  & \textbullet &  &  & \textbullet &  &  & Analytical &  \\
Salameh \& Jaber \cite{salameh2000economic} & 2000 & Imperfect & EOQ &  &  & \textbullet &  & \textbullet &  &  &  & Analytical &  \\
Tripathy et al. \cite{tripathy2003eoq} & 2003 & Imperfect & EOQ &  &  & \textbullet &  & \textbullet &  &  &  & Analytical &  \\
Pasandideh et al. \cite{pasandideh2011genetic} & 2011 &  & SC & \textbullet &  & \textbullet &  &  & \textbullet & \textbullet &  &  & GA \\
Pasandideh et al. \cite{pasandideh2011parameter} & 2011 &  & SC & \textbullet &  &  & \textbullet &  & \textbullet & \textbullet &  &  & GA \\
Taleizadeh et al. \cite{taleizadeh2011multiple} & 2011 &  & SC &  & \textbullet &  & \textbullet &  & \textbullet & \textbullet & \textbullet &  & HS,GA \\
Chen et al. \cite{chen2014retailer} & 2014 &  & SC & \textbullet &  & \textbullet &  & \textbullet &  &  &  & Analytical &  \\
Pasandideh et al. \cite{pasandideh2014optimization} & 2014 &  & SC & \textbullet &  & \textbullet &  &  & \textbullet & \textbullet &  & Solver & GA \\
Cárdenas-Barrón \& Sana \cite{cardenas2015multi} & 2017 &  & SC & \textbullet &  & \textbullet &  & \textbullet &  &  &  & Analytical &  \\
Khan et al. \cite{khan2017learning} & 2018 & Imperfect & SC & \textbullet &  & \textbullet &  & \textbullet &  &  &  & Analytical &  \\
Mokhtari \cite{mokhtari2018joint} & 2018 & Reusable & EOQ &  &  & \textbullet &  & \textbullet &  &  &  & Analytical &  \\
Tiwari et al. \cite{tiwari2018joint} & 2018 & Deteriorating & SC & \textbullet &  & \textbullet &  & \textbullet &  & \textbullet &  & Analytical &  \\
Karimian et al. \cite{karimian2020geometric} & 2020 &  & SC & \textbullet &  & \textbullet &  &  & \textbullet & \textbullet & \textbullet & \begin{tabular}[c]{@{}l@{}}Geometric\\ Programming\end{tabular} &  \\
Pourmohammad-Zia et al. \cite{pourmohammad2021food} & 2021 & Growing & SC & \textbullet &  &  & \textbullet & \textbullet &  &  &  & Game theory &  \\
Pourmohammad-Zia et al. \cite{pourmohammad2021dynamic} & 2021 & Growing and deteriorating & SC & \textbullet &  & \textbullet &  & \textbullet &  &  &  & Analytical &  \\
Mokhtari et al. \cite{mokhtari2022economic} & 2022 & Substitutable & EPQ &  &  & \textbullet &  &  & \textbullet & \textbullet &  & Solver &  \\
Asadkhani et al. \cite{asadkhani2022sustainable} & 2022 & Imperfect & SC & \textbullet &  & \textbullet &  & \textbullet &  &  &  & Analytical &  \\
Fallahi et al. \cite{fallahi2022constrained} & 2022 & Reusable & EOQ &  &  & \textbullet &  &  & \textbullet & \textbullet &  & Interior point & \begin{tabular}[c]{@{}l@{}}DE,PSO,\\ DEQL,PSOQL\end{tabular} \\
\hline
This article & 2023 & Reusable & SC & \textbullet &  &  & \textbullet &  & \textbullet & \textbullet & \textbullet & \begin{tabular}[c]{@{}l@{}}Sequential Quadratic\\ Programming\end{tabular}  & GWO,WOA \\
\hline 
\multicolumn{14}{l}{\begin{tabular}[l]{@{}l@{}}GA: Genetic algorithm, HS: Harmony Search, DE: Differential Evolution, PSO: Particle Swarm Optimization, DEQL: Differential Evolution-Q-Learning, PSOQL: Particle Swarm Optimization-Q Learning,\\GWO: Grey Wolf Optimizer, WOA: Whale Optimization Algorithm\end{tabular}}\\
\end{tabular}
\end{adjustbox}
\end{table}

To the best of our knowledge, no other research focuses on the inventory management of reusable items in a two-level supply chain under operational constraints. The goal of this research is to address this problem and propose a new problem that helps the supply chain managers of reusable items to coordinate the vendor with the retailer through the determination of optimal inventory decisions for the integrated systems. This problem is presented as a single-vendor multi-product multi-retailer inventory system under operational constraints. Also, several sources of uncertainty may impact the constraints of inventory systems. In this paper, we assume that the system's constraints are stochastic and handle it by the chance-constrained programming method. Additionally, GWO and WOA novel metaheuristic algorithms are designed and implemented as the solution approach. The efficiency and effectiveness of these algorithms are shown by validating the results using the SQP algorithm as a powerful exact method.

\section{Problem presentation and mathematical modeling} \label{sec:problem}
In this section, we will present the new problem for inventory management of reusable products in a single-vendor multi-product multi-retailer supply chain and formulate the mathematical model.

\subsection{Notations}
Let us consider the following notations:

\begin{table}[H]
\begin{adjustbox}{scale = 0.95}
\centering
\begin{tabular}{p{1.5cm} p{15cm}}
\hline
Sets & \\
\hline \hline
$J$ & Retailer index; $j \in \{1,\dots,J\}$ \\
$K$ & Item index; $k \in \{1,\dots,K\}$ \\
\hline 
Parameters & \\
\hline \hline
$ILU_{jk}(t)$ & The inventory level of usable item $k$ of retailer $j$ \\
$ILR_{jk}(t)$ & The inventory level of recoverable item $k$ of retailer $j$ \\
$OCS_{jk}$ & The ordering cost of vendor per order of reusable item $k$ from retailer $j$ \\
$OCU_{jk}$ & The ordering cost of retailer $j$  per order of reusable item $k$ \\
$OCR_{jk}$ & The fixed recovery cost of retailer $j$ per recovery of recoverable item $k$ \\
$PC_{k}$ & The unit purchasing cost of reusable item $k$ with mean  $\mu_k^{PC}$ and standard deviation $\sigma_k^{PC}$ \\
$RC_{jk}$ & The unit recovery cost of recoverable item $k$ for retailer $j$ \\
$HCU_{jk}$ & The unit holding cost of usable item $k$ per unit of time for retailer $j$ with mean $\mu_{jk}^{HCU}$ and standard deviation $\sigma_{jk}^{HCU}$ \\
$HCR_{jk}$ & The unit holding cost of recoverable item $k$ per unit of time for retailer $j$ with mean $\mu_{jk}^{HCR}$  and standard deviation $\sigma_{jk}^{HCR}$ \\
$D_{jk}$ & The demand rate of reusable item $k$ for retailer $j$ with mean $\mu_{jk}^D$  and standard deviation $\sigma_{jk}^D$ \\
$m_k$ & The maximum number that reusable item $k$ can be reused and recovered \\
$f_k$ & The required storage capacity for storing reusable item $k$ with mean 
$\mu_k^f$ and standard deviation $\sigma_k^f$ \\
$B_j$ & The maximum budget of retailer $j$ with mean $\mu_j^B$ and standard deviation $\sigma_j^B$  \\ 
$AHU_j$	& The maximum holding cost for usable items of retailer $j$ with mean $\mu_j^{AHU}$ and standard deviation $\sigma_j^{AHU}$ \\
$AHR_j$	& The maximum holding cost for recoverable items of retailer $j$ with mean $\mu_j^{AHR}$ and standard deviation $\sigma_j^{AHR}$ \\
$WSU_j$ & The maximum storage capacity of retailer $j$ for usable items with mean $\mu_j^{WSU}$ and standard deviation $\sigma_j^{WSU}$ \\
$WSR_j$ & The maximum storage capacity of retailer $j$ for recoverable items with mean $\mu_j^{WSR}$ and standard deviation $\sigma_j^{WSR}$ \\
$WS$ & The total maximum storage capacity of the vendor with mean $\mu^{WS}$ and standard deviation $\sigma^{WS}$ \\
$N$ & The maximum number of orders for all items with mean $\mu^N$ and standard deviation $\sigma^N$ \\
$\alpha$ & The probability of violating each of the constraints \\
\hline 
Variables & \\
\hline \hline
$Q_{jk}$ & The economic order quantity for reusable item $k$ of retailer $j$ per cycle \\
$q_{jk}$ & The economic reuse and recovery quantity for reusable item $k$ of retailer $j$ per cycle \\
$p_{jk}$ & The ratio of economic order quantity to economic reuse and recovery quantity for reusable item $k$ of retailer $j$ per cycle \\
$TCB_j$ & The total cost of retailer $j$ \\
$TCS$ & The total cost of the vendor \\
$TCE$ & The total cost of the supply chain \\
\hline
\end{tabular}
\end{adjustbox}
\end{table}

\subsection{Assumptions}
The main assumption of the presented problem can be expressed as follows:
\begin{itemize}
    \item There are one vendor, $K$ item, and $J$ retailer in the system.
    \item The demand rate for products is constated and deterministic.
	\item The maximum storage capacity of usable products for each retailer is less than an upper limit with a probability greater than or equal to $\alpha$.
	\item The maximum storage capacity of recoverable products for each retailer is less than an upper limit with a probability greater than or equal to $\alpha$.
	\item The maximum holding cost of usable products for each retailer is less than an upper limit with a probability greater than or equal to $\alpha$.
	\item The maximum holding cost of recoverable products for each retailer is less than an upper limit with a probability greater than or equal to $\alpha$.
	\item The maximum budget for each retailer is less than an upper limit with a probability greater than or equal to $\alpha$.
	\item The total number of orders in the system is less than an upper limit with a probability greater than or equal to $\alpha$.
	\item There is no lead time in the system.
	\item Backorder and lost sale shortages are not allowed. 
\end{itemize}

\subsection{Problem definition}
\begin{figure}[H]
\centering 
\includegraphics[width=17cm]{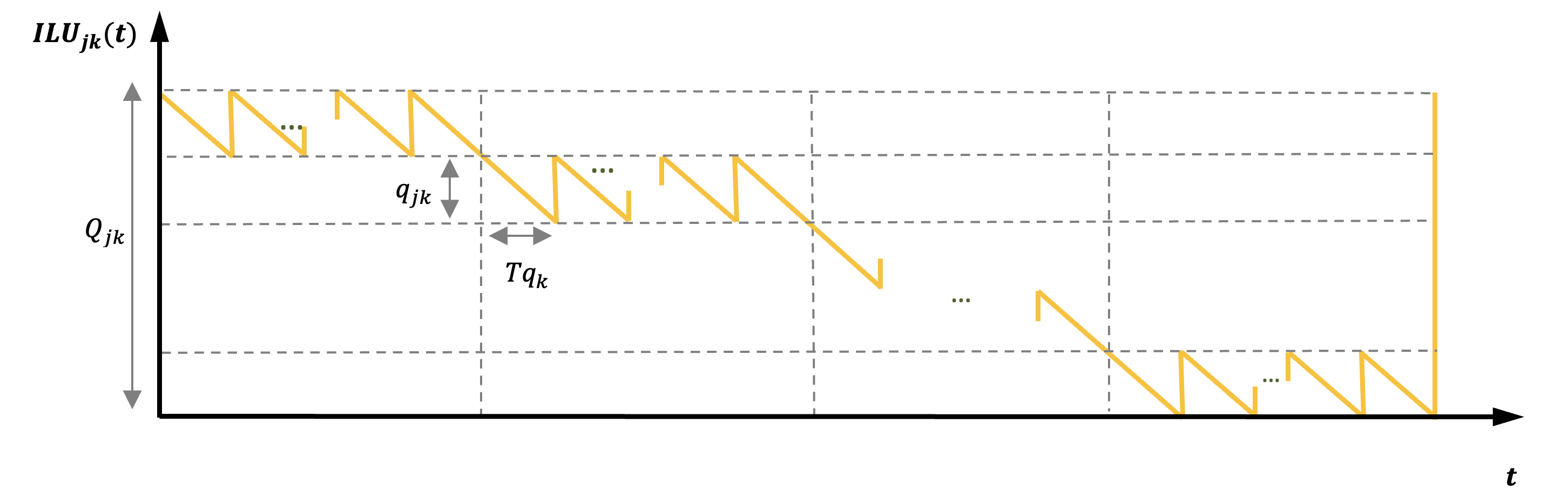}
\caption{The stock level of usable product $k$ in the warehouse of retailer $j$}
\label{fig:ilu}
\end{figure}

\begin{figure}[H]
\centering 
\includegraphics[width=17cm]{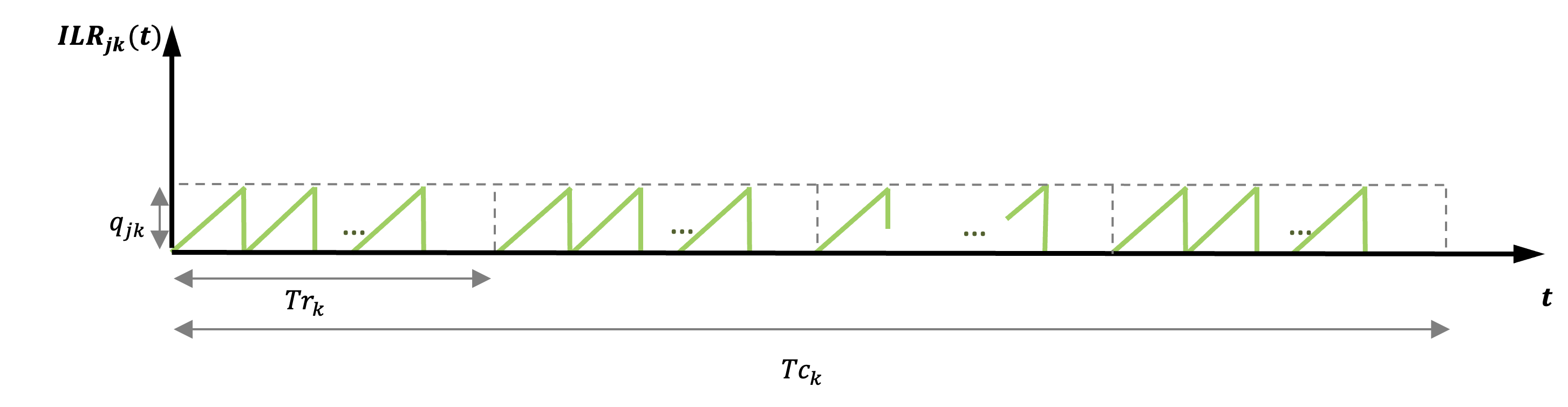}
\caption{The istock level of recoverable product $k$ in the warehouse of retailer $j$}
\label{fig:ilr}
\end{figure}

Consider a two-level supply chain of reusable items, including a vendor and $j \in \{1,\dots,J\}$ retailers. In this system, the usable term refers to the products ready to satisfy customers' demands. In addition, recoverable products are the products that need a recovery process to become usable for demand satisfaction. Each retailer needs to place orders for $k \in \{1,\dots,K\}$ reusable items from the vendor. The retailer $j$ purchases the reusable item $k$ from the vendor at the unit purchasing cost $PC_k$. For each order, $OCU_{jk}$ is the ordering cost of retailer $j$ for reusable item $k$. In addition, $OCS_{jk}$ is the imposed ordering cost to the vendor regarding the order of retailer $j$ for reusable item $k$. Retailer $j$ uses the purchased reusable item $k$ to satisfy the demand of customers $D_{jk}$. As mentioned before, the products are reusable, and the used products can be recovered and used again for a maximum of $m_k$ times. The recovery cost of each unit of product $k$ for retailer $j$ is $RC_{jk}$. In addition, there is a fixed recovery operational cost for the recovery of product $k$ by retailer $j$, which is $OCR_{jk}$. The presence of usable and recoverable product $k$ in warehouses of retailer $j$ imposes holding costs on each retailer, which are specified by $HCU_{jk}$ and $HCR_{jk}$, respectively. The stock level diagrams of usable and recoverable product $k$ in the warehouses of retailer $j$ are shown in Figures~\ref{fig:ilu} and \ref{fig:ilr}.

A set of operational constraints are considered in the system to bring it into the real-world environment. The constraints are stochastic, and it is assumed that the resources are available with a probability equal to and greater than $\alpha$. The limitation on the maximum storage capacity of usable and recoverable products, the holding costs of usable and recoverable products, the total available purchasing of usable and recoverable products, and the total number of orders, are the stochastic constraints of the supply chain. 

\subsection{Mathematical modeling}
Regarding the explained problem, the cost components of the supply chain can be described as follows:
\begin{itemize}
    \item The total purchasing cost of retailer $j$:
    \begin{equation}
        \sum_{k=1}^{K} PC_k \frac{D_{jk}}{(m_k+1)} \mkern53mu \quad \forall j \in \mathcal{J}
    \end{equation}
    \item The total fixed ordering cost of the vendor:
    \begin{equation}
        \sum_{j=1}^{J}\sum_{k=1}^{K} OCS_{jk} \frac{D_{jk}}{(m_k+1)p_{jk}q_{jk}}
    \end{equation}
    \item The total fixed ordering cost of retailer $j$:
    \begin{equation}
        \sum_{k=1}^{K} OCU_{jk} \frac{D_{jk}}{(m_k+1)p_{jk}q_{jk}} \quad \forall j \in \mathcal{J}
    \end{equation}
    \item The total fixed recovery cost of retailer $j$:
    \begin{equation}
        \sum_{k=1}^{K} OCR_{jk} D_{jk} (\frac{m_k}{m_k+1}) \quad \forall j \in \mathcal{J}
    \end{equation}
    \item The total recovery operational cost of retailer $j$:
    \begin{equation}
        \sum_{k=1}^{K} RC_{jk} \frac{D_{jk}}{q_{jk}} (\frac{m_k}{m_k+1}) \quad \forall j \in \mathcal{J}
    \end{equation}
    \item The total holding cost of usable products of retailer $j$:
    \begin{equation}
        \sum_{k=1}^{K} HCU_{jk} (\frac{p_{jk}q_{jk}}{2}) \quad \forall j \in \mathcal{J}
    \end{equation}
    \item The total holding cost of recoverable products of retailer $j$:
    \begin{equation}
        \sum_{k=1}^{K} HCR_{jk} (\frac{m_k}{m_k+1}) \frac{q_{jk}}{2} \quad \forall j \in \mathcal{J}
    \end{equation}
\end{itemize}

Considering the above-described components, the total cost objective of retailer $j$ can be expressed as follows:
\begin{align}
& \nonumber
TCB_j(p_{jk}, q_{jk}) = \sum_{k=1}^{K} PC_k \frac{D_{jk}}{(m_k+1)} + \sum_{k=1}^{K} OCU_{jk} \frac{D_{jk}}{(m_k+1)p_{jk}q_{jk}} \\&  \nonumber
+\sum_{k=1}^{K} OCR_{jk} D_{jk} (\frac{m_k}{m_k+1}) + \sum_{k=1}^{K} RC_{jk} \frac{D_{jk}}{q_{jk}} (\frac{m_k}{m_k+1}) \\& 
+\sum_{k=1}^{K} HCU_{jk} (\frac{p_{jk}q_{jk}}{2}) + \sum_{k=1}^{K} HCR_{jk} (\frac{m_k}{m_k+1}) \frac{q_{jk}}{2} & 
\end{align}

Also, the vendor bears the following cost:
\begin{align}
&  TCS(p_{jk}, q_{jk}) = \sum_{j=1}^{J}\sum_{k=1}^{K} OCS_{jk} \frac{D_{jk}}{(m_k+1)p_{jk}q_{jk}} &
\end{align}

In the decentralized systems, each retailer places the order regarding her cost function. Here, the total cost of other retailers and the vendor do not play any role in the determination of the replenishment decisions. As a consequence of such a policy, a huge cost may impact the overall performance of the system. In centralized decision-making, managers try to determine the optimal decisions regarding the costs of all entities. In this situation, a coordination mechanism should be used to establish the integrated total cost function of the supply chain. The presented mechanism by Hill \cite{hill1997single} is one of the well-known mechanisms managers tend to adapt to coordinate such supply chains \cite{asadkhani2022sustainable,moheb2023reverse}. Therefore, we use this coordination approach to integrate the cost components of the vendor and retailer. The centralized total cost of the two-layer network under the coordination mechanism of Hill \cite{hill1997single} can be expressed as follows:
\begin{align}
& \nonumber
TCE = \sum_{j=1}^{J}\sum_{k=1}^{K} OCS_{jk} \frac{D_{jk}}{(m_k+1)p_{jk}q_{jk}} + \sum_{k=1}^{K} PC_k \frac{D_{jk}}{(m_k+1)} \\&  \nonumber
+ \sum_{k=1}^{K} OCU_{jk} \frac{D_{jk}}{(m_k+1)p_{jk}q_{jk}} +\sum_{k=1}^{K} OCR_{jk} D_{jk} (\frac{m_k}{m_k+1}) \\& \nonumber
+ \sum_{k=1}^{K} RC_{jk} \frac{D_{jk}}{q_{jk}} (\frac{m_k}{m_k+1}) +\sum_{k=1}^{K} HCU_{jk} (\frac{p_{jk}q_{jk}}{2})\\& 
+ \sum_{k=1}^{K} HCR_{jk} (\frac{m_k}{m_k+1}) \frac{q_{jk}}{2} & 
\end{align}

The integrated objective function is subjected to the following stochastic operational constraints:

\begin{align}
    & Min TCE  \label{objective}\\
    & \text{s.t.} \notag\\
    & {P\left(\sum_{k=1}^K PC_k.p_{jk}.q_{jk} \leq B_j\right) \geq 1-\alpha} \quad && \forall j\in \mathcal{J} \label{con1}\\
    & {P\left(\sum_{j=1}^J\sum_{k=1}^K f_k.p_{jk}.q_{jk} \leq WS\right) \geq 1-\alpha} \label{con2}\\
    & {P\left(\sum_{k=1}^K f_k.p_{jk}.q_{jk} \leq WSU_j\right) \geq 1-\alpha} \quad && \forall j\in \mathcal{J} \label{con3}\\
    & {P\left(\sum_{k=1}^K f_k.q_{jk} \leq WSR_j\right) \geq 1-\alpha} \quad && \forall j\in \mathcal{J} \label{con4}\\
    & {P\left(\sum_{k=1}^K HCU_{jk}.(\frac{p_{jk}.q_{jk}}{2}) \leq AHU_j  \right) \geq 1-\alpha} \quad && \forall j\in \mathcal{J} \label{con5}\\
    & {P\left(\sum_{k=1}^K HCR_{jk}.(\frac{m_k}{m_k+1}).\frac{q_{jk}}{2} \leq AHR_j  \right) \geq 1-\alpha} \quad && \forall j\in \mathcal{J} \label{con6}\\
    & {P\left(\sum_{j=1}^J\sum_{k=1}^K \frac{D_{jk}}{(m_k+1).p_{jk}.q_{jk}} \leq N  \right) \geq 1-\alpha} \label{con7}\\
    & {p_{jk}, q_{jk} \geq 0 } \quad && \forall j\in \mathcal{J} , k \in \mathcal{K} \label{con8}
\end{align}

Constraints \eqref{con1} to \eqref{con7} are the chance constraints of the system. Stochastic constraints \eqref{con1} ensure that the total purchasing cost of the products for each retailer does not exceed the retailer's maximum budget. Constraint \eqref{con2} limits the total available storage capacity for the vendor. Stochastic constraints \eqref{con3} and \eqref{con4} specify the storage capacity of each retailer for usable and recoverable products. The constraints of the maximum holding cost of usable and recoverable for each retailer is shown via stochastic constraints \eqref{con5} and \eqref{con6}. Stochastic constraints \eqref{con7} limit the system's total number of orders. Constraints \eqref{con8} determine the type of decision variables. 

We utilize the chance constraint programming approach to deal with the uncertainty of the constraints. Considering a normal probability distribution with mean $\mu$ and standard deviation $\sigma$ for the upper bound of each stochastic constraint, the constraints can be rewritten as follows \cite{pasandideh2015optimization}:
\begin{align}
    & Min TCE  \label{objective2}\\
    & \text{s.t.} \notag\\
    & {\sum_{k=1}^K \mu_k^{PC}.p_{jk}.q_{jk} + Z_\alpha .\sqrt{\sum_{k=1}^K (\sigma_k^{PC}.p_{jk}.q_{jk})^2 +(\sigma_j^B)^2} \leq \mu_j^B} \quad && \forall j\in \mathcal{J} \label{conR1}\\
    & {\sum_{j=1}^J \sum_{k=1}^K \mu_k^f.p_{jk}.q_{jk} +Z_\alpha .\sqrt{\sum_{j=1}^J \sum_{k=1}^K (\sigma_k^f.p_{jk}.q_{jk})^2 +(\sigma^{WS})^2} \leq \mu^{WS}} \label{conR2}\\
    & {\sum_{k=1}^K \mu_k^f.p_{jk}.q_{jk} + Z_\alpha .\sqrt{\sum_{k=1}^K (\sigma_k^f.p_{jk}.q_{jk})^2 +(\sigma_j^{WSU})^2} \leq \mu_j^{WSU}} \quad && \forall j\in \mathcal{J} \label{conR3}\\
    & {\sum_{k=1}^K\mu_k^f.q_{jk} +Z_\alpha.\sqrt{\sum_{k=1}^K(\sigma_k^f.q_{jk})^2 +(\sigma_j^{WSR})^2} \leq \mu_j^{WSR}} \quad && \forall j\in \mathcal{J} \label{conR4}\\
    & {\sum_{k=1}^K\mu_{jk}^{HCU}.(\frac{p_{jk}.q_{jk}}{2}) +Z_\alpha.\sqrt{\sum_{k=1}^K \left( \sigma_{jk}^{HCU}.(\frac{p_{jk}.q_{jk}}{2}) \right)^2 +(\sigma_j^{AHU})^2} \leq \mu_j^{AHU}} \quad && \forall j\in \mathcal{J} \label{conR5}\\
    & {\sum_{k=1}^K\mu_{jk}^{HCR}.(\frac{m_k}{m_k+1}).\frac{q_{jk}}{2} + Z_\alpha.\sqrt{\sum_{k=1}^K \left( \sigma_{jk}^{HCR}.(\frac{m_k}{m_k+1}).\frac{q_{jk}}{2} \right)^2 +(\sigma_j^{AHR})^2} \leq \mu_j^{AHR}} \quad && \forall j\in \mathcal{J} \label{conR6}\\
    & {\sum_{j=1}^J\sum_{k=1}^K \frac{\mu_{jk}^D}{(m_k+1).p_{jk}.q_{jk}} +Z_\alpha.\sqrt{\sum_{j=1}^J\sum_{k=1}^K \left( \frac{\sigma_{jk}^D)}{(m_k+1).p_{jk}.q_{jk}} \right)^2 +(\sigma^N)^2}  \leq \mu^N } \label{conR7} \\
    & {p_{jk}, q_{jk} \geq 0 } \quad \forall j\in \mathcal{J} , \forall k \in \mathcal{K}
\end{align}

In the above equations, the upper $\alpha$-percentile point of the normal probability distribution (standard form) is shown by $Z_\alpha$.

The presented model is a constrained nonlinear programming model, which can not easily be solved via exact classical methods or commercial solvers. Therefore, we utilize the metaheuristic algorithms as the solution methodology. The algorithms will be presented in the next section.

\section{Solution approach} \label{sec:sol}
This section will discuss the proposed solution methodologies to solve the problem. The presented model is a constrained nonlinear programming mathematical model. Previous research pointed out that such inventory management problems are challenging to solve with classical methods due to the nonlinearity of the model and several local optimum solutions \cite{fallahi2022constrained,khalilpourazari2016optimization}. Therefore, metaheuristic algorithms are widely used as a powerful solution for multi-product inventory management in supply chains \cite{fallahi2022constrained,pasandideh2011genetic,taleizadeh2011multiple}. This paper uses GWO and WAO as two recently developed metaheuristics to solve the problem. In addition, we use SQP as an exact approach to show the efficiency of the metaheuristics.

\subsection{Grey wolf optimizer metaheuristic algorithm}
GWO is a nature-inspired population-based metaheuristic that was introduced by Mirjalili et al. \cite{mirjalili2014grey}, and extensively used as the solution approach to optimization problems in several fields \cite{faris2018grey,jiang2022dsgwo}. The GWO metaheuristic is designed to find the near global optimum of a given function for the solution spaces with continuous variables. The algorithm is inspired from the hunting behavior of grey wolves in a pack, where each wolf plays a specific role in cooperating and competing with each other to find the best prey. GWO is initialized with a population of random solutions, and then it iteratively updates the solutions by mimicking the interactions among the wolves in the pack. The steps of the GWO metaheuristic can be summarized are as below:

\subsubsection{Initialize prey and hunters}
The first step of the GWO metaheuristic is to initialize the population of solutions $X=\{x_1,x_2,\dots,x_n\}$ where $x_i$ is a solution vector in the search space. The population is initialized with a set of random solutions. The solutions number, also known as the population size, is a user-defined parameter that can vary depending on the problem at hand. Larger population size can increase the diversity of solutions but also increases the computational cost. 
GWO also starts with initializing the position and fitness of the alpha, beta, and delta wolves, denoted by $x_\alpha$, $x_\beta$ and $x_\delta$, respectively. These wolves are used as reference points in the next steps of the algorithm. The $x_\alpha$ wolf is considered the leader of the pack and has the best objective function in the population. The $x_\beta$ wolf is the second-best solution, and the $x_\delta$ wolf is the solution with the third-best objective function. After the initialization of the population in the first iteration, there is a set of common steps in the next iterations of GWO, which is the core of the algorithm. These steps iteration process are repeated until a stopping criterion is met, such as the maximum number of iterations or achieving a satisfactory solution. In each iteration, the algorithm updates the position of the $x_\alpha$, $x_\beta$, and $x_\delta$ wolves, as well as the position of the other solutions in the population, based on the details that are as follows.

\subsubsection{Hunting}
\begin{figure}[H]
\centering 
\includegraphics[width=17cm]{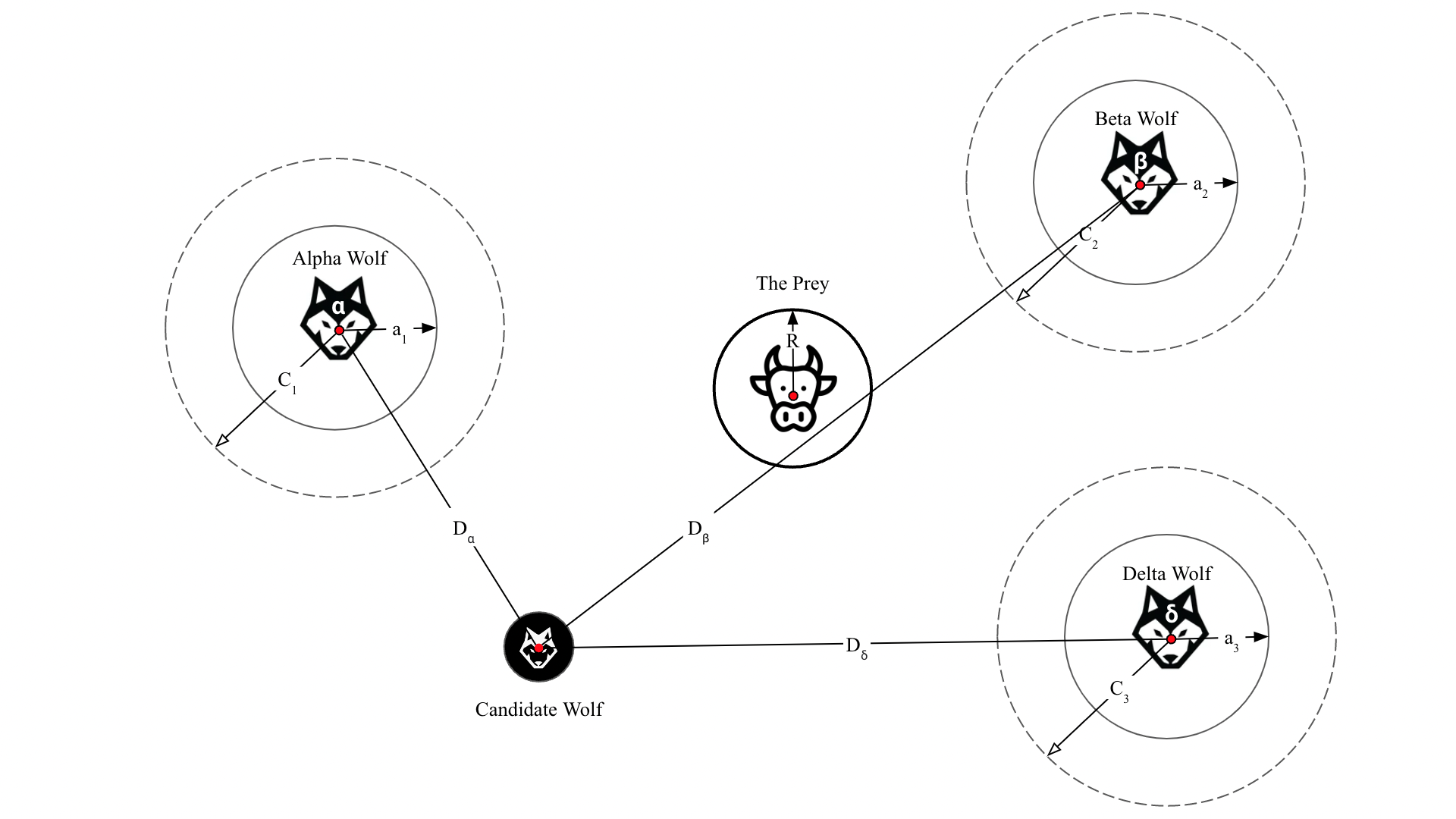}
\caption{The search process of the GWO metaheuristic algorithm}
\label{fig:gwo}
\end{figure}

GWO updates the position of the solutions based on the hunting behavior of grey wolves, where the wolves cooperate and compete to find the best prey. The position of the $x_\alpha$, $x_\beta$, and $x_\delta$ wolves is updated using the following equations:
\begin{equation} \label{eq:update-alpha-wolf}
    x_{t+1}^\alpha= x_t^\alpha+\alpha_t (x_t^\beta-x_t^\alpha )+\beta_t (x_t^\delta-x_t^\alpha )
\end{equation}
\begin{equation} \label{eq:update-beta-wolf}
    x_{t+1}^\beta= x_t^\beta+\alpha_t (x_t^\alpha-x_t^\beta )+\beta_t (x_t^\delta-x_t^\beta )
\end{equation}
\begin{equation} \label{eq:update-delta-wolf}
    x_{t+1}^\delta= x_t^\delta+\alpha_t (x_t^\alpha-x_t^\delta )+\beta_t (x_t^\beta-x_t^\delta )
\end{equation}

where $\alpha_t$ and $\beta_t$ are linearly decreasing functions of the iteration $t$, and are used to control the step size of the search. In addition, the position of the other solutions in the population is updated using the following equation:
\begin{equation}
    x_{t+1}^{(i)}= x_t^{(i)}+r.(x_{t+1}^\alpha-x_t^{(i)} ) + r.(x_{t+1}^\beta-x_t^{(i)})+r.(x_{t+1}^\delta-x_t^{(i)} )
\end{equation}

where $r$ is a generator that generates a random number between 0 and 1. After updating the position of the solutions, the fitness of the $x_\alpha$, $x_\beta$, and $x_\delta$ wolves is re-evaluated, and the new $x_\alpha$, $x_\beta$, and $x_\delta$ wolves are selected from the population. The search process of the algorithm is shown in Figure~\ref{fig:gwo}.

\subsubsection{Stop Criteria}
The stopping criterion is a pre-defined condition that is used to determine when the search process should stop. The stopping criterion is usually based on the number of iterations or the quality of the solutions. A common stopping criterion is to stop the algorithm after a certain number of iterations, also known as the maximum number of iterations. This criterion is used to prevent the algorithm from running indefinitely, and it is typically set based on the computational resources available and the complexity of the problem. Another stopping criterion is to stop the algorithm when a satisfactory solution is found. This criterion is used to stop the optimization process when the solution quality reaches a certain level. We set the first criterion as the stop criterion of GWO. Algorithm~\ref{alg:gwo} presents the pseudocode of the GWO.

\begin{algorithm}[hbt!]
\caption{GWO Algorithm}\label{alg:gwo}
\begin{algorithmic}[1]
\STATE \textbf{Input:} Maximum iteration ($t_{max}$), Population size ($N_{pop}$), $a$, $A$, and $C$
\STATE \textbf{Output:} Best solution ($\Vec{X}_\alpha$)
\FOR{$i= 1:N_{pop}$}
\STATE Initialize GWO solution $\Vec{X}_i^0$
\STATE Calculate the fitness $f(\Vec{X}_i^0)$
\ENDFOR
\STATE $\Vec{X}_\alpha \gets $ The first best search agent
\STATE $\Vec{X}_\beta \gets $ The second best search agent
\STATE $\Vec{X}_\delta \gets $ The third best search agent
\FOR{$t= 1:t_{max}$}
\STATE $\Vec{A} \gets 2\Vec{a}.\Vec{r}_1 - \Vec{a}$
\STATE $\Vec{C} \gets 2\Vec{r}_2$
\STATE $\Vec{D}_\alpha \gets |\Vec{C}_1.\Vec{X}_\alpha-\Vec{X}|,\Vec{D}_\beta \gets |\Vec{C}_2.\Vec{X}_\beta-\Vec{X}|,\Vec{D}_\delta \gets |\Vec{C}_3.\Vec{X}_\delta-\Vec{X}|$
\STATE $\Vec{X}_1^t \gets \Vec{X}_\alpha-\Vec{A}_1.(\Vec{D}_\alpha ),\Vec{X}_2^t \gets \Vec{X}_\beta-\Vec{A}_2.(\Vec{D}_\beta ),\Vec{X}_3^t \gets \Vec{X}_\delta-\Vec{A}_3.(\Vec{D}_\delta )$
\STATE $\Vec{A} \gets 2.\Vec{a}.\Vec{r}_1 - \Vec{a}$
\STATE $\Vec{X}^{t+1} \gets \frac{\Vec{X}_1^t+\Vec{X}_2^t+\Vec{X}_3^t}{3}$ \COMMENT{Update the position of the search agent}
\STATE Update $a$, $A$, and $C$ multipliers
\STATE Calculate the fitness $f(\Vec{X}^{t+1})$
\STATE Update $x_\alpha$, $x_\beta$, and $x_\delta$ using equations \eqref{eq:update-alpha-wolf} - \eqref{eq:update-delta-wolf}
\STATE $t \gets t+1$
\ENDFOR
\STATE \textbf{Return:} $\Vec{X}_\alpha$ 
\end{algorithmic}
\end{algorithm}

\subsection{Whale Optimization Algorithm}
WOA is another nature-inspired population-based metaheuristic optimization algorithm inspired by the foraging behavior of humpback whales. The algorithm was proposed by Mirjalili \& Lewis \cite{mirjalili2016whale}, and it is used to solve optimization problems in several research areas \cite{zhao2022multipopulation,chakraborty2021enhanced}. It is known for its ability to find high-quality solutions, as well as its ability to avoid getting stuck in local optima \cite{gharehchopogh2019comprehensive}. As pointed out, WOA is a population-based optimization algorithm, where a group of candidate solutions, called a population, is iteratively improved to converge towards a near-optimal solution. The algorithm simulates the foraging behavior of humpback whales, where each whale represents a candidate solution, and the search space is divided into subproblems. The main steps of WOA are:

\subsubsection{Initialize prey}
The first step of the WOA algorithm is to randomly initialize a population of candidate solutions. The population is typically a set of $n$ solutions, where each solution is represented by a vector of $d$ variables, $x_i^T=\{x_i^{(1)},x_i^{(2)},\dots,x_i^{(d)}\}$. The values of the variables should be chosen within the specified bounds. 

\subsubsection{Hunting}
\begin{figure}[H]
\centering 
\includegraphics[width=17cm]{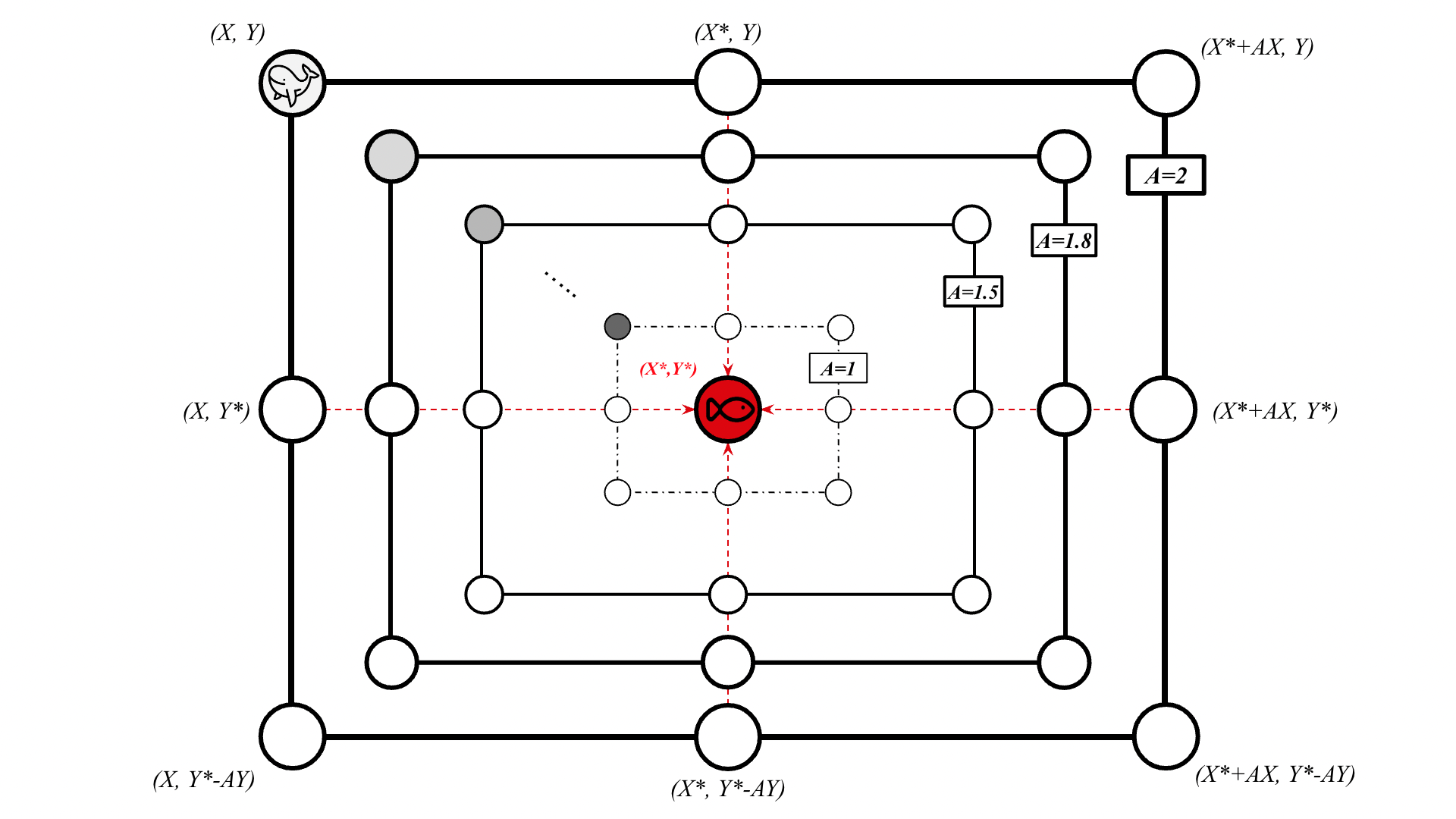}
\caption{The shrinking updating process of the WOA metaheuristic algorithm}
\label{fig:woa1}
\end{figure}

In the WOA, the current global best solution, also known as the leader, is assumed to be close to the optimum, and the other solutions in the population are guided towards it, similar to how humpback whales encircle their prey. In other words, the leader solution is used as a point of reference for the other solutions to follow, guiding the search of the other solutions toward better regions of the space. The other solutions are updated simultaneously, based on their distance and fitness difference from the leader solution. This behavior is represented by the mathematical equations~\eqref{eq:update-whale1} and~\eqref{eq:update-whale2} that are used to update the positions of the solutions in the population. These equations are designed to mimic the foraging behavior of humpback whales.
\begin{equation} \label{eq:update-whale1}
    \Vec{D} = |\Vec{C}.\Vec{X}_t^* - \Vec{X}_t|
\end{equation}
\begin{equation} \label{eq:update-whale2}
    \Vec{X}_{t+1} = \Vec{X}_t^* - \Vec{A}.\Vec{D}
\end{equation}

where $\Vec{X}_t^*$ is the best solution obtained in $t^{th}$ iteration of WOA. The parameters $\Vec{A}$ and $\Vec{C}$ are as follows:
\begin{equation} \label{eq:whale-A}
    \Vec{A} = 2.\Vec{a}.\Vec{r}_1 - \Vec{a}
\end{equation}
\begin{equation} \label{eq:whale-C}
    \Vec{C} = 2.\Vec{r}
\end{equation}

where $\Vec{a}$ decreases linearly from 2 to 0 over the iterations. Also, $\Vec{r}$ is a random number generated uniformly between 0 and 1. Two models are used to represent the attacking behavior of humpback whales: 

\paragraph{Shrinking updating position}
This behavior is accomplished by decreasing the value of $\Vec{a}$ in equation~\eqref{eq:whale-A}. Note that the range of fluctuation is also decreased by $\Vec{A}$. Figure~\ref{fig:woa1} illustrates the potential positions that can be reached from $(X,Y)$ to $(X^*,Y^*)$ when $0 \leq A \leq 1$ in a 2-dimensional space.

\begin{figure}[H]
\centering 
\includegraphics[width=5cm]{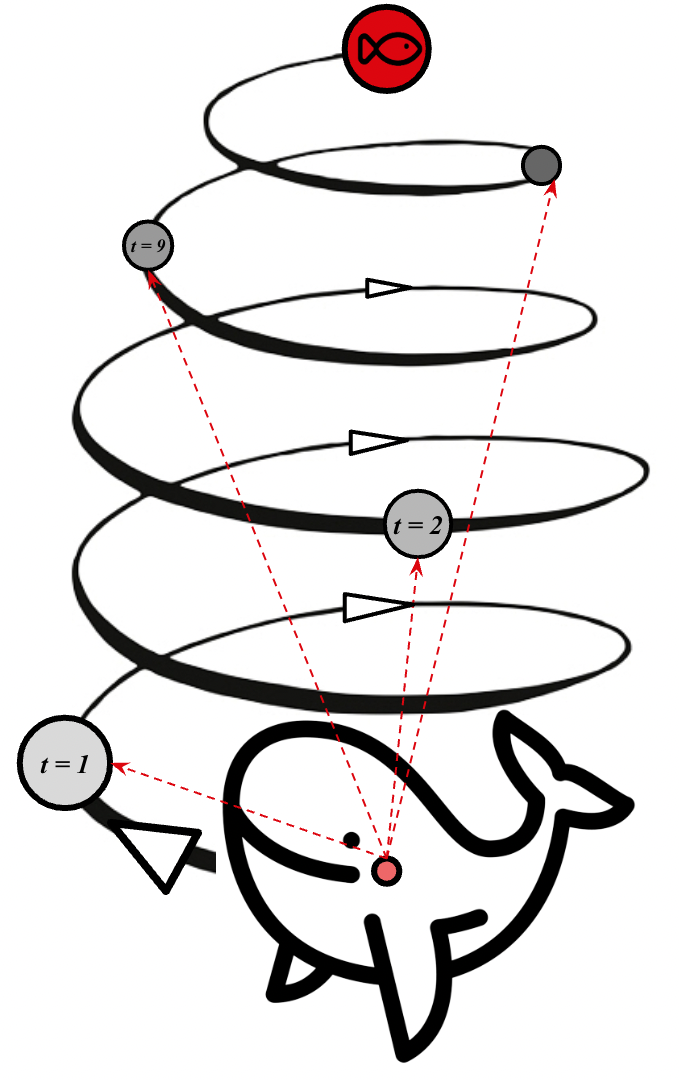}
\caption{The spiral updating process of the WOA metaheuristic algorithm}
\label{fig:woa2}
\end{figure}

\paragraph{Spiral updating position}
Another observation of Humpback whales' hunting is swimming in a helical path toward their prey \cite{mirjalili2016whale}. To replicate this behavior, a spiral function is defined to modify the position of search as:
\begin{equation} \label{eq:whale-Spiral}
    \Vec{X}_{t+1} = \Vec{D}^\prime . e^{bt} . \cos{(2 \pi l)} + \Vec{X}_t^*
\end{equation}
where the distance of the $i^{th}$ whale to the prey (best solution found so far) is represented by $|\Vec{X}_t^* - \Vec{X}_t|$, $b$ is a constant used to shape the logarithmic spiral, $0 \leq l \leq 1$ is a random number. The procedure is shown in Figure~\ref{fig:woa2}.

To simulate how humpback whales move around their prey by swimming in a shrinking circle and along a spiral-shaped path, the algorithm utilizes a probability of 50\% to either use the first or second method.

\subsubsection{Stop criteria}
WOA is to repeat the process of evaluating the fitness, selecting the leader, and updating the positions of the solutions until a stopping criterion is met. We set the maximum number of iterations as the stopping criterion. Algorithm~\ref{alg:woa} shows the pseudocode of the WOA.

\begin{algorithm}[hbt!]
\caption{WOA Algorithm}\label{alg:woa}
\begin{algorithmic}[1]
\STATE \textbf{Input:} Maximum iteration ($t_{max}$), Population size ($N_{pop}$), $a$, $b$, $A$, $C$, $p$, and $l$
\STATE \textbf{Output:} Best solution ($\Vec{X}^*$)
\FOR{$i= 1:N_{pop}$}
\STATE Initialize whale position $\Vec{X}_i^0$
\STATE Calculate the fitness $f(\Vec{X}_i^0)$
\ENDFOR
\FOR{$i= 1:t_{max}$}
\FOR{$Each search agent$}
\STATE Update $a$, $A$, $C$, $p$, and $l$
\IF{$p < 0.5 $}
\STATE $\Vec{A}  \gets  2.\Vec{a}.\Vec{r}_1 - \Vec{a}$
\STATE $\Vec{C}  \gets  2.\Vec{r}$
\IF{$|\Vec{A}| < 1$}
\STATE $\Vec{D} \gets |\Vec{C}.\Vec{X}_t^* - \Vec{X}_t|$
\STATE $\Vec{X}_{t+1}  \gets  \Vec{X}_t^* - \Vec{A}.\Vec{D}$
\ELSIF{$|\Vec{A}| \geq 1$}
\STATE Select a random search agent ($\Vec{X}_{rand}$)
\STATE $\Vec{D} \gets |\Vec{C}.\Vec{X}_{rand} - \Vec{X}_t|$
\STATE $\Vec{X}_{t+1}  \gets  \Vec{X}_{rand} - \Vec{A}.\Vec{D}$
\ENDIF
\ELSIF{$p \geq 0.5 $}
\STATE $\Vec{D}^\prime \gets |\Vec{X}^* - \Vec{X}|$
\STATE $\Vec{X}_{t+1}  \gets  \Vec{D}^\prime . e^{bt} . \cos{(2 \pi l)} + \Vec{X}_t^*$
\ENDIF
\ENDFOR
\STATE Check and amend if any search agent goes beyond the search space
\STATE Calculate the fitness $f(\Vec{X}^{t+1})$
\STATE Update $\Vec{X}_t^*$ if there is a better solution
\STATE $t \gets t+1$
\ENDFOR
\STATE \textbf{Return:} $\Vec{X}^*$ 
\end{algorithmic}
\end{algorithm}

\subsection{Sequential quadratic programming}
The SQP algorithm is a technique for solving nonlinear optimization problems involving smooth and nonsmooth functions. It is an iterative method that uses a combination of gradient and Hessian information to determine the next iterate. SQP is particularly well-suited for solving large-scale NLP problems and has been shown to be effective in many applications. 

The algorithm is based on the theory of Quadratic Programming (QP), and it's a combination of gradient and Hessian information to determine the next iterate \cite{gill2011sequential,gurwitz1989sequential}. It uses the Karush-Kuhn-Tucker (KKT) conditions to manage equality constraints in the same way that Newton's technique does when solving an unconstrained NLP optimization problem. \cite{rockafellar1970convex}. The KKT conditions are a set of necessary and sufficient conditions that a solution to a constrained optimization problem must satisfy. The solution of QP sub-problem is typically used to determine a line search direction in the SQP algorithm. SQP is similar to the active-set algorithm and has some advantages over other exact solution methods. One advantage is that the SQP method guarantees exact feasibility with respect to bounds. This means that the algorithm will always find a feasible solution that satisfies all the bounds constraints. Another advantage of SQP is that it is more robust to problems with complex values \cite{rockafellar1970convex, bartlett2000active}. This is because the SQP algorithm approximates the objective function and constraints, which can help avoid getting stuck in poor local solutions and help the algorithm converge to a global optimum.

SQP is also used in the literature to determine the economic order (production) quantity in constrained multi-product inventory problems \cite{pasandideh2015optimization,khalilpourazari2016optimization,khalilpourazari2019modeling}. The SQP method is particularly well-suited to solving this type of problem because it can handle the nonlinear and nonconvex nature of the objective function and constraints that arise in the EPQ problem with stochastic constraints. In addition, studies have shown that the SQP method can perform significantly better than other approaches, such as the interior-point method. This is because the SQP method can often converge to a global optimum, whereas the interior-point method can get stuck in poor local solutions. Also, the SQP method can handle the nonlinear and nonconvex nature of the objective function and constraints that arise in the EPQ problem with stochastic constraints more effectively than the interior-point exact approach \cite{pasandideh2015optimization,khalilpourazari2016optimization}.

\section{Computational experiments}\label{sec:result}
In this section, we will evaluate the model's and metaheuristic algorithms' performance by solving the numerical examples. We use the data from the works by Mokhtari \cite{mokhtari2018joint}, and Fallahi et al. \cite{fallahi2022constrained} to generate numerical examples. The details of the data are presented in Table ~\ref{tab:datex}. 
\begin{table}[H]
\caption{The data of numerical examples}
\label{tab:datex}
\centering
\begin{tabular}{cccc}
\hline
Parameter & Range & Parameter & Range \\ \hline
\textit{OCS} & \textit{~UN(1300,1900)} & \textit{f} & \textit{~UN(1,2)} \\
\textit{OCU} & \textit{~UN(1300,1900)} & \textit{B} & \textit{~UN(290000000,310000000)} \\
\textit{OCR} & \textit{~UN(80,100)} & \textit{AHU} & \textit{~UN(380000,420000)} \\
\textit{PC} & \textit{~UN(40,60)} & \textit{AHR} & \textit{~UN(1900000,2100000)} \\
\textit{RC} & \textit{~UN(16,24)} & \textit{WSR} & \textit{~UN(18000,22000)} \\
\textit{HCU} & \textit{~UN(1,2)} & \textit{WSU} & \textit{~UN(18000,22000)} \\
\textit{HCR} & \textit{~UN(6,10)} & \textit{WS} & \textit{25000*J} \\
\textit{D} & \textit{~UN(10000,14000)} & \textit{N} & \textit{10000*K} \\
\textit{m} & \textit{~DUN(2,5)} &  &  \\ \hline
\end{tabular}
\end{table}

The algorithms are run on a personal laptop with 16 GB Ram and an Intel Core i7 4.7 GHz CPU. We also provide the SAS code for SQP solver.

\begin{figure}
\centering 
\includegraphics[width=10cm]{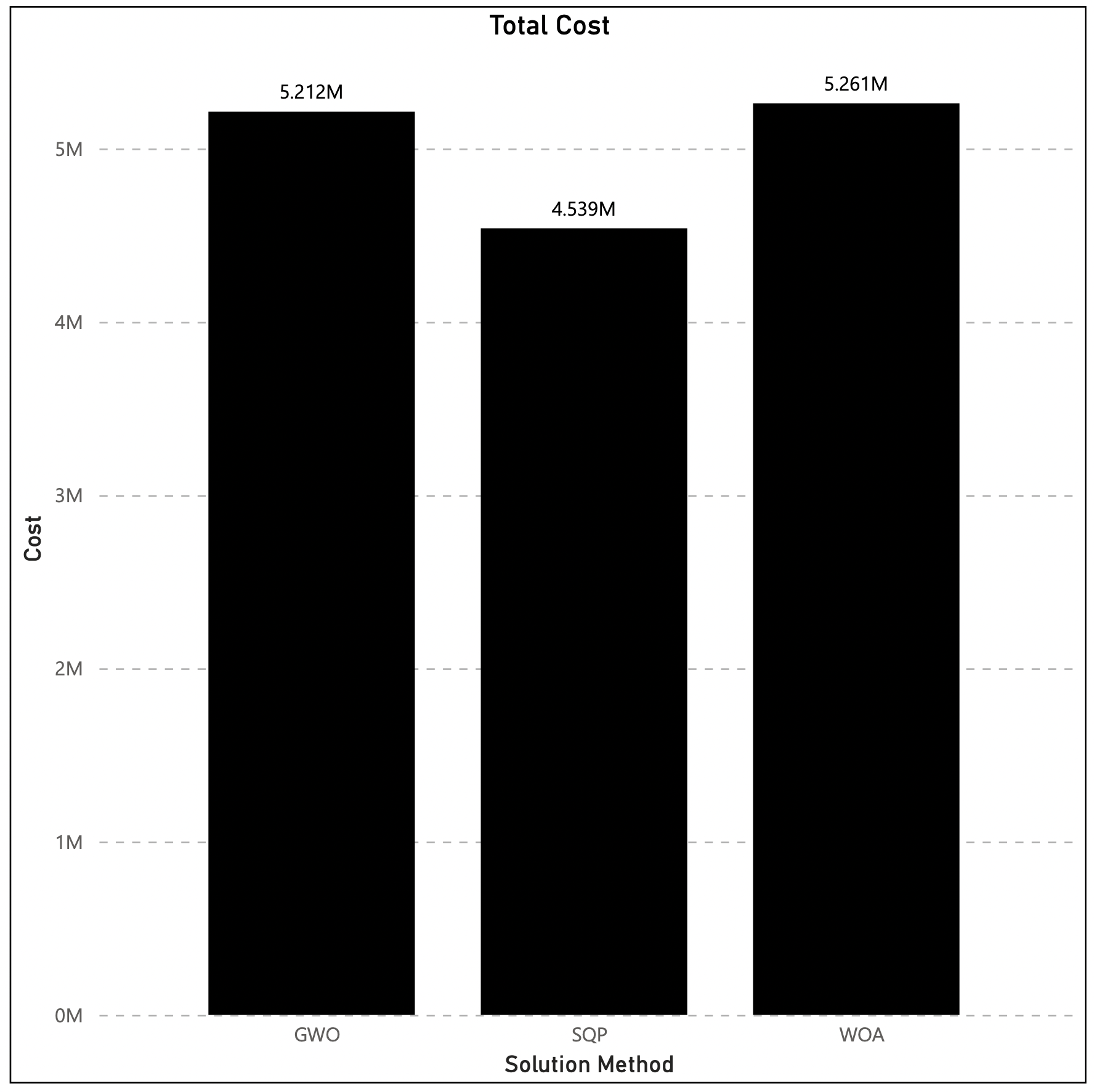}
\caption{Comparison between the value of objective functions using different solution methods.}
\label{fig:c-total}
\end{figure}

First, we are going to validate the performance of metaheuristic algorithms by comparing the results of a small-size numerical example. For this goal, we consider a numerical example with two retailers and one product in the system.  Figure~\ref{fig:c-total} shows the calculated results of algorithms for the numerical example. As obvious, there is no significant difference between the performance of the metaheuristics. The total costs of GWO and WOA are more than SQP by about 14.78 and 15.90, respectively. Such difference confirms that the algorithms perform well in searching the solution space. As obvious, GWO has a better performance than WOA.

\begin{figure*}[hbt!]
    \centering
    \begin{subfigure}[b]{0.3\textwidth}
        \centering
        \includegraphics[width=\textwidth]{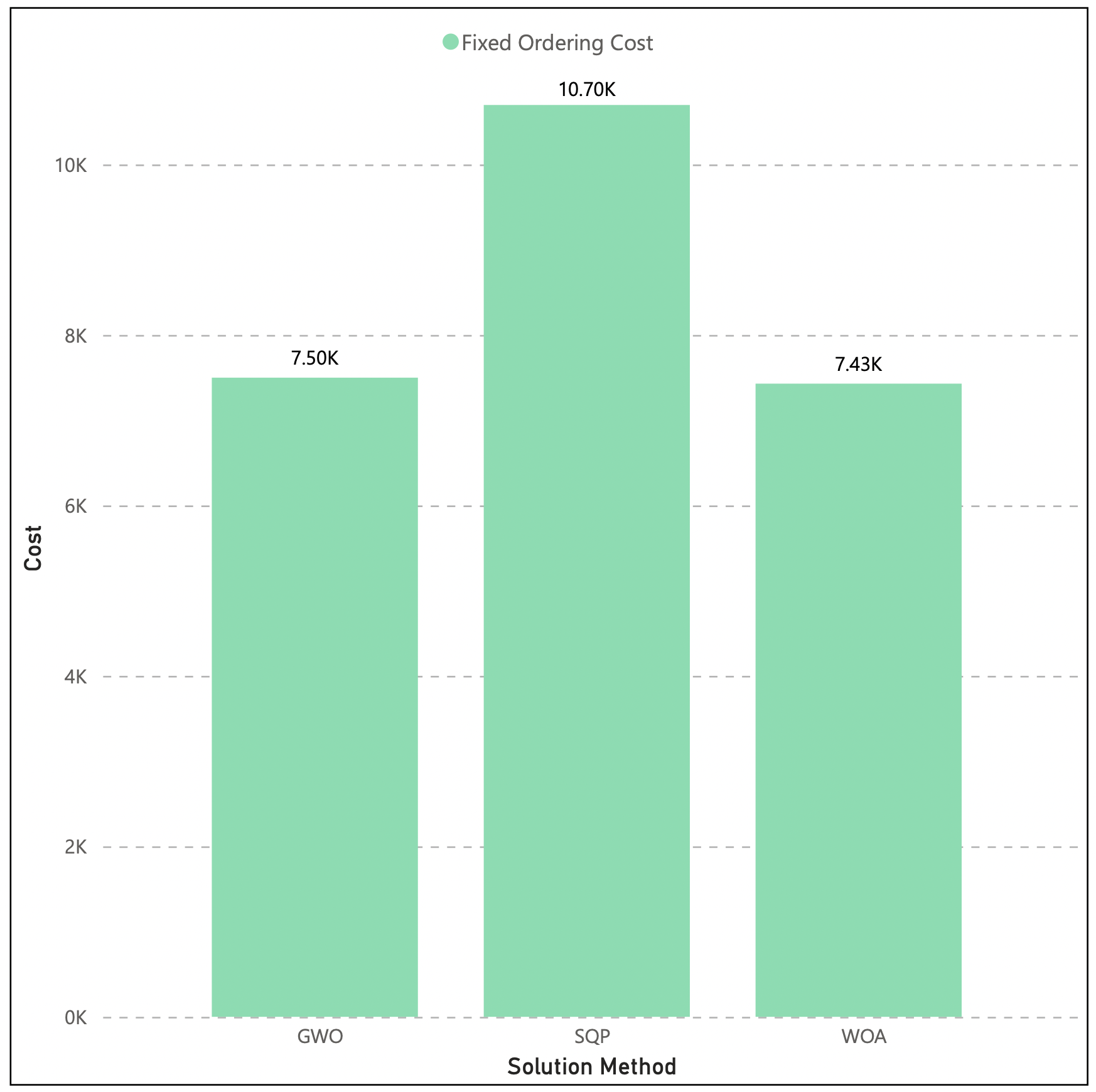}
        \caption[]
        {{\small Fixed Ordering Cost}}    
        \label{fig:c1}
    \end{subfigure}
    \hfill
    \begin{subfigure}[b]{0.3\textwidth}  
        \centering 
        \includegraphics[width=\textwidth]{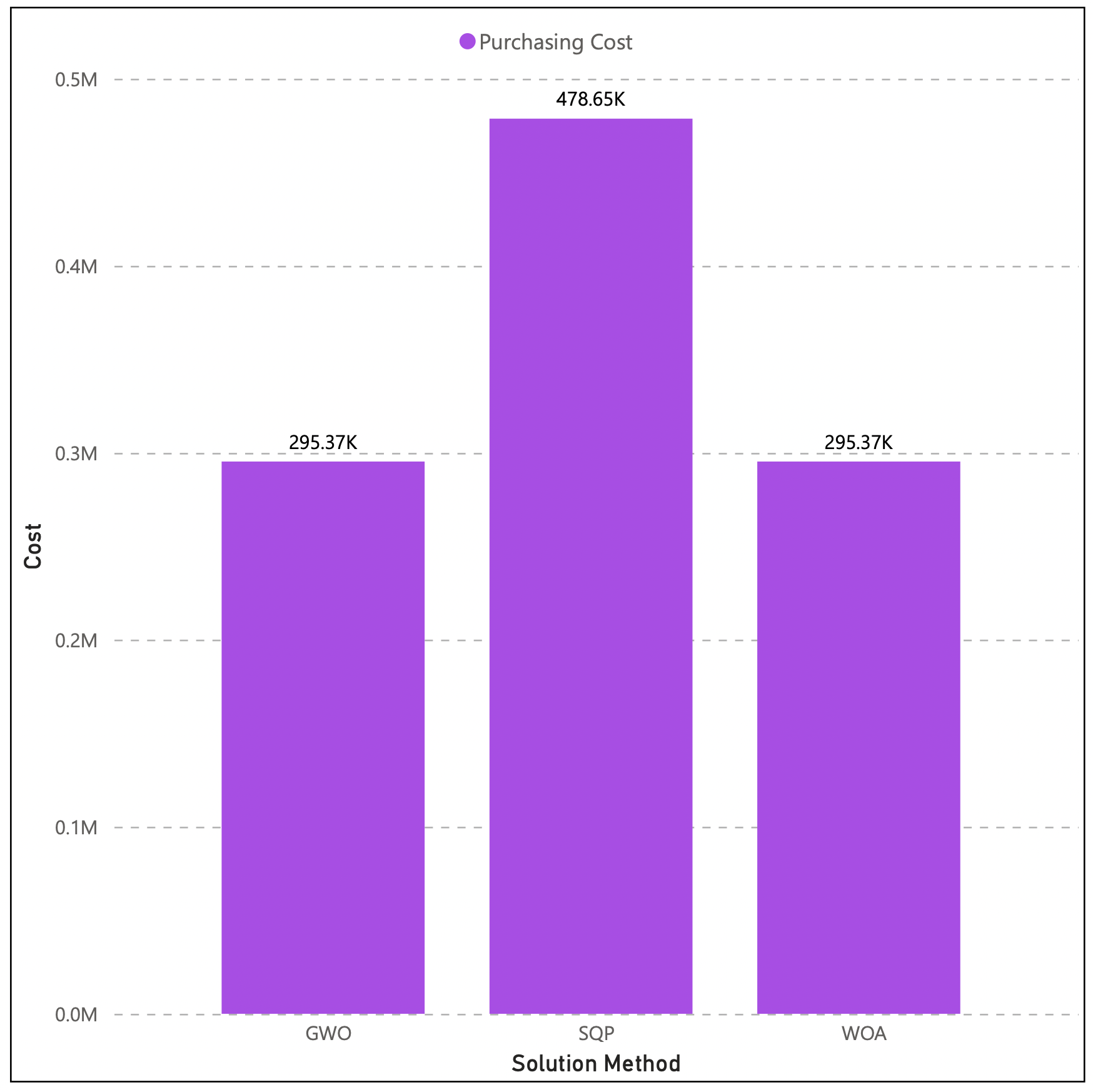}
        \caption[]
        {{\small Purchasing Cost}}    
        \label{fig:c2}
    \end{subfigure}
    \hfill
    \begin{subfigure}[b]{0.3\textwidth}  
        \centering 
        \includegraphics[width=\textwidth]{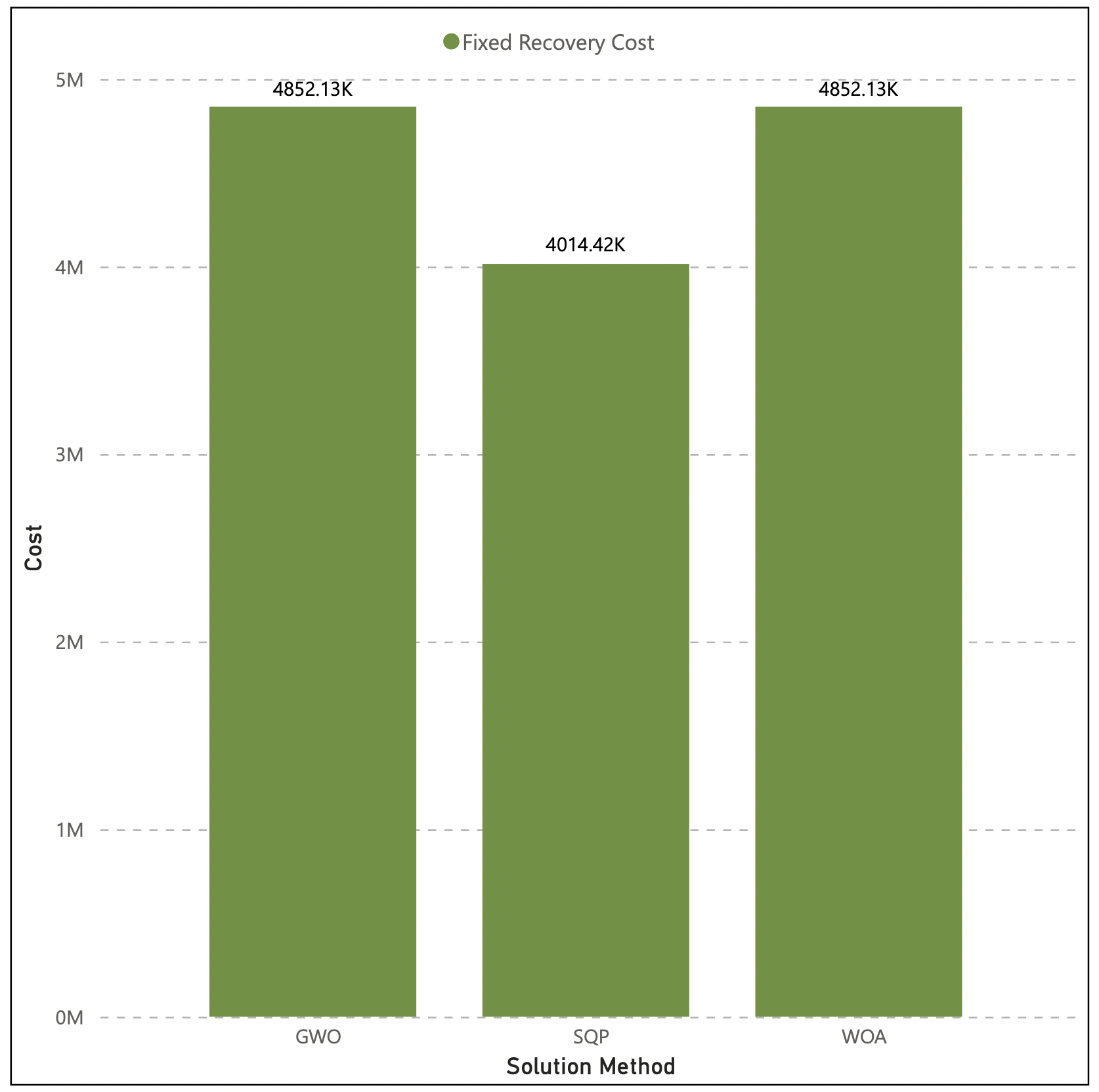}
        \caption[]
        {{\small Recovery Cost (Fixed)}}    
        \label{fig:c3}
    \end{subfigure}    
    \vskip\baselineskip
    \begin{subfigure}[b]{0.3\textwidth}   
        \centering 
        \includegraphics[width=\textwidth]{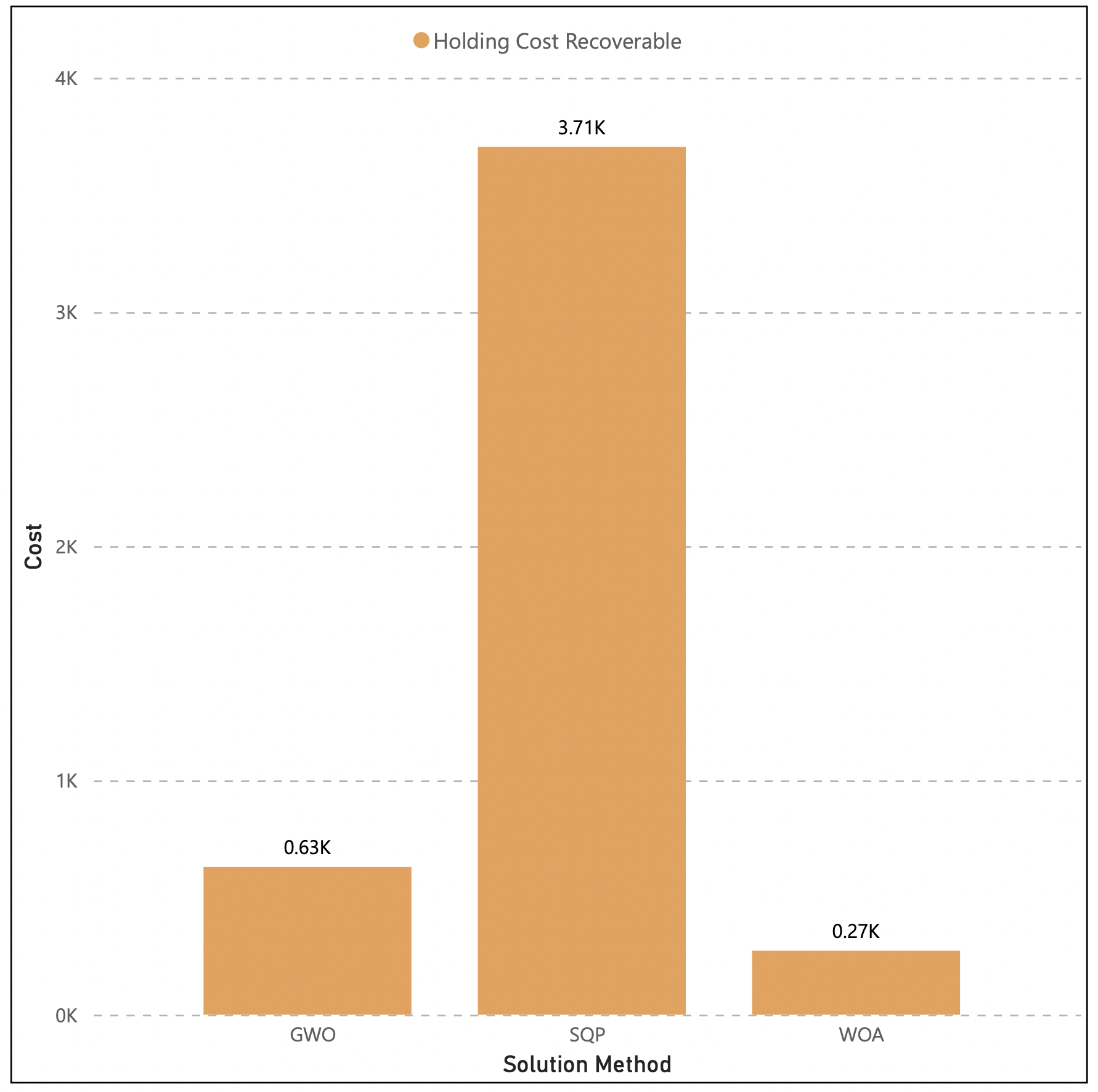}
        \caption[]
        {{\small Holding Cost (Recoverable)}}    
        \label{fig:c4}
    \end{subfigure}
    \hfill
    \begin{subfigure}[b]{0.3\textwidth}   
        \centering 
        \includegraphics[width=\textwidth]{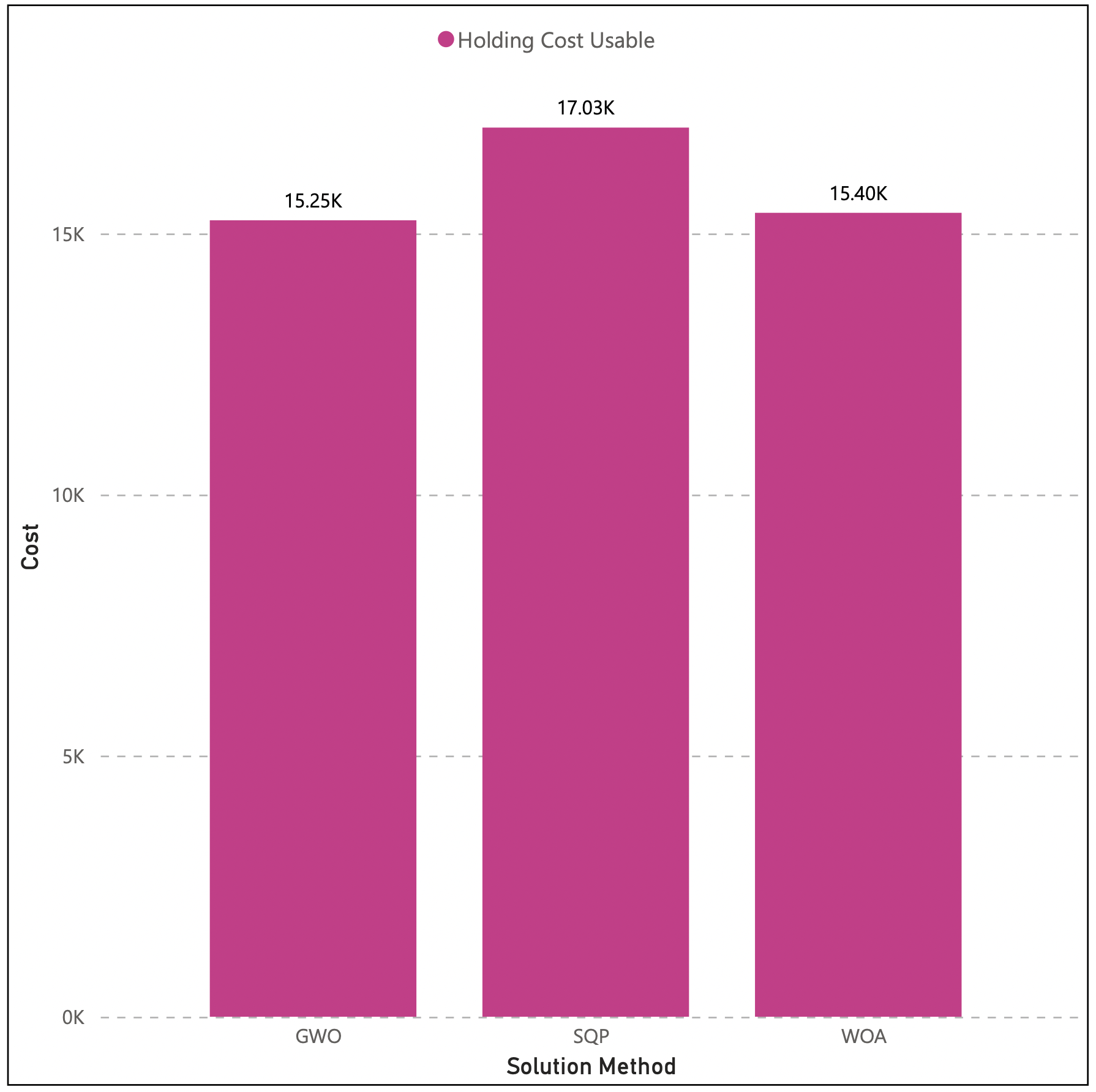}
        \caption[]
        {{\small Holding Cost (Usable)}}    
        \label{fig:c5}
    \end{subfigure}
    \hfill
    \begin{subfigure}[b]{0.3\textwidth}   
        \centering 
        \includegraphics[width=\textwidth]{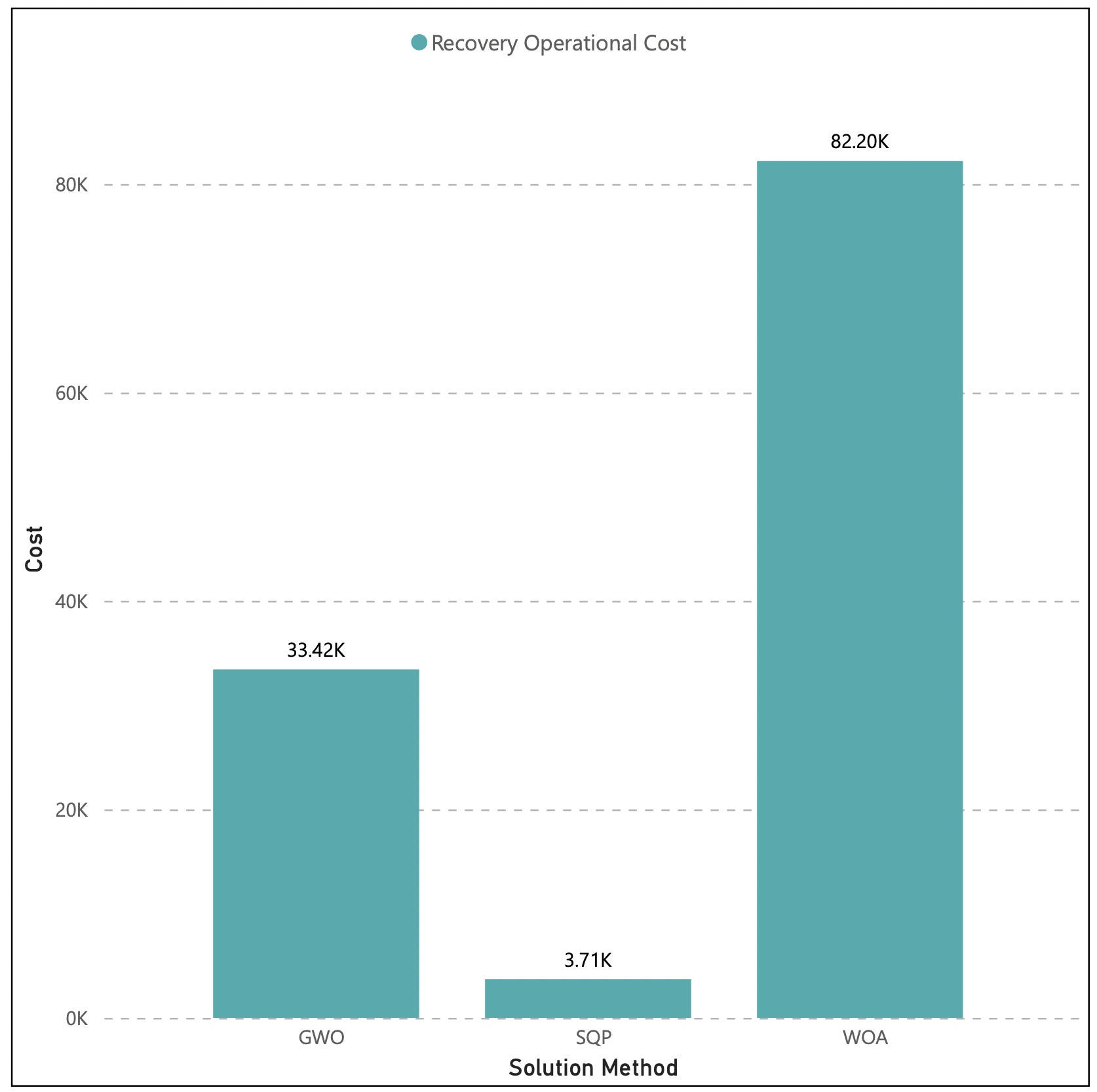}
        \caption[]
        {{\small Recovery Cost (Operational)}}    
        \label{fig:c6}
    \end{subfigure}    
    \caption[]
    {\small The calculated cost components of GWO, WOA, and SQP algorithms} 
    \label{fig:cost-terms}
\end{figure*}

To provide better insight, the cost component by each algorithm is also provided in Figure~\ref{fig:cost-terms}. As can be seen, a great portion of the total cost is due to the fixed recovery cost. In addition, the holding cost of recoverable products is less than the other cost components of the system.

\subsection{Parameter tuning}
The input parameters highly impact the performance of metaheuristic algorithms \cite{mousavi2021intelligent}. Various methods are employed in the literature to determine the input parameters of metaheuristics. The trial-and-error methods are very time-consuming and do not guarantee the quality of solutions. Therefore, using a systematic approach for parameter tuning seems necessary. Taguchi's design \cite{roy2010primer} of experiments is one of the widely used statistical methods for this goal. The Taguchi method utilizes the concept of the orthogonal array to manage the number of experiments. In this statistical approach, the affecting factors are grouped into two categories of signal ($S$) and noise ($N$) factors. There is no direct control over the noise factors, and they can not easily be changed or removed. Therefore, Taguchi tries to find the optimum level of signal factors in such as way that the effect of noise factors is minimized. Taguchi defines the relative importance of individual components in terms of their primary effects on the objective function in order to determine the best parameter levels. The repeated data is transformed by Taguchi into a different value, which is the variation measure. This transformation is signal-to-noise (S⁄N) ratio, which is calculated as below for a minimization problem:
\begin{equation} \label{eq:taguchi}
    S⁄N=-10 \log{\frac{1}{m} \sum_{j=1}^m z_j^2}
\end{equation}
where $m$ is the total replications and $z_j$ is the response in $j^{th}$ replications. Three considered levels for each parameter of algorithms are presented in Table~\ref{tab:taguchi}. The L9 orthogonal arrays are utilized for the parameter calibration of GWO and WOA. Three levels are defined for each parameter based on the details in Table~\ref{tab:taguchi}. In addition, each orthogonal array is run in five replications.

\begin{table}[]
\caption{The considered levels for parameters of metaheuristic algorithms}
\label{tab:taguchi}
\centering
\begin{tabular}{llcccccc}
\hline
\multicolumn{1}{c}{\multirow{3}{*}{Parameter}} &                      & \multicolumn{2}{c}{GWO} &  & \multicolumn{3}{c}{WOA} \\ \cline{3-4} \cline{6-8} 
\multicolumn{1}{c}{}                           &                      & A           & B         &  & A       & B      & C    \\ \cline{3-4} \cline{6-8} 
\multicolumn{1}{c}{}                           &                      & Max\_it     & N\_pop    &  & Max\_it & N\_pop & b    \\ \hline
Level 1                                        & \multicolumn{1}{c}{} & 100         & 100       &  & 100     & 100    & -0.9 \\
Level 2                                        & \multicolumn{1}{c}{} & 150         & 150       &  & 150     & 150    & -1   \\
Level 3                                        & \multicolumn{1}{c}{} & 200         & 200       &  & 200     & 200    & -1.1 \\ \hline
\end{tabular}
\end{table}

The optimal level of parameters are presented in Table~\ref{tab:taguchi-optimum}. 

\begin{table}[H]
\caption{The optimal input parameters of metaheuristic algorithms}
\label{tab:taguchi-optimum}
\centering
\begin{tabular}{llcccccc}
\hline
\multicolumn{1}{c}{\multirow{3}{*}{Parameter}} &                      & \multicolumn{2}{c}{GWO} &  & \multicolumn{3}{c}{WOA} \\ \cline{3-4} \cline{6-8} 
\multicolumn{1}{c}{}                           &                      & A           & B         &  & A       & B      & C    \\ \cline{3-4} \cline{6-8} 
\multicolumn{1}{c}{}                           &                      & Max\_it     & N\_pop    &  & Max\_it & N\_pop & b    \\ \hline
Optimal level                                        & \multicolumn{1}{c}{} & 200         & 200       &  & 200     & 200    & -1.1 \\ \hline
\end{tabular}
\end{table}

\begin{figure}[H]
\centering 
\includegraphics[width=17cm]{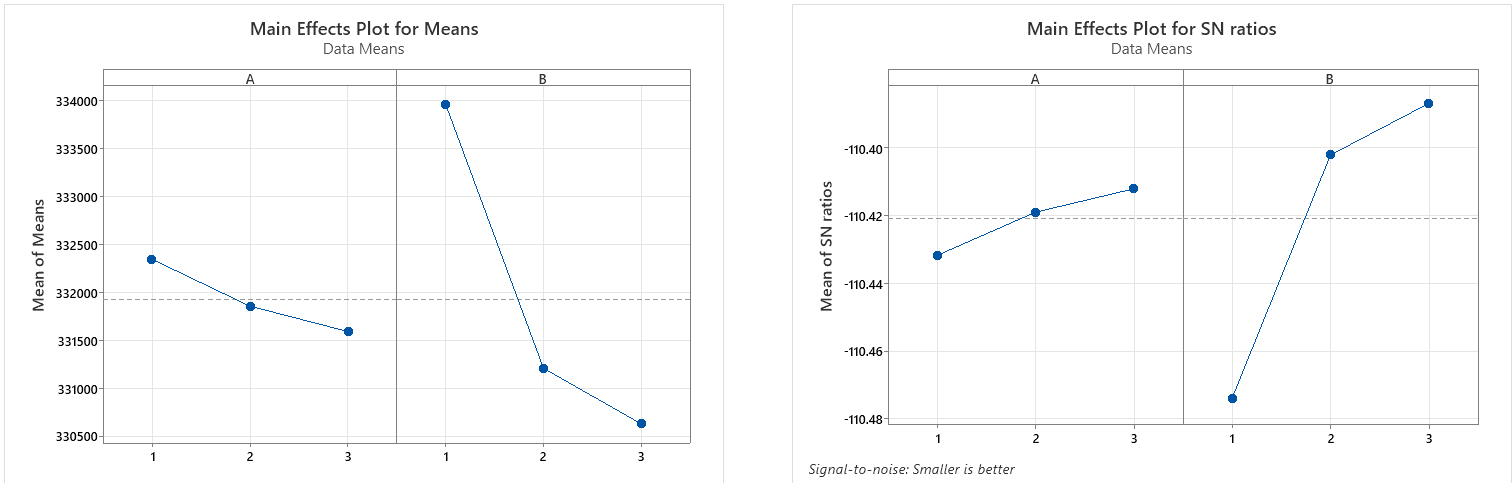}
\caption{The main effects plot for S⁄N of GWO.}
\label{fig:taguchi-gwo}
\end{figure}

\begin{figure}[H]
\centering 
\includegraphics[width=17cm]{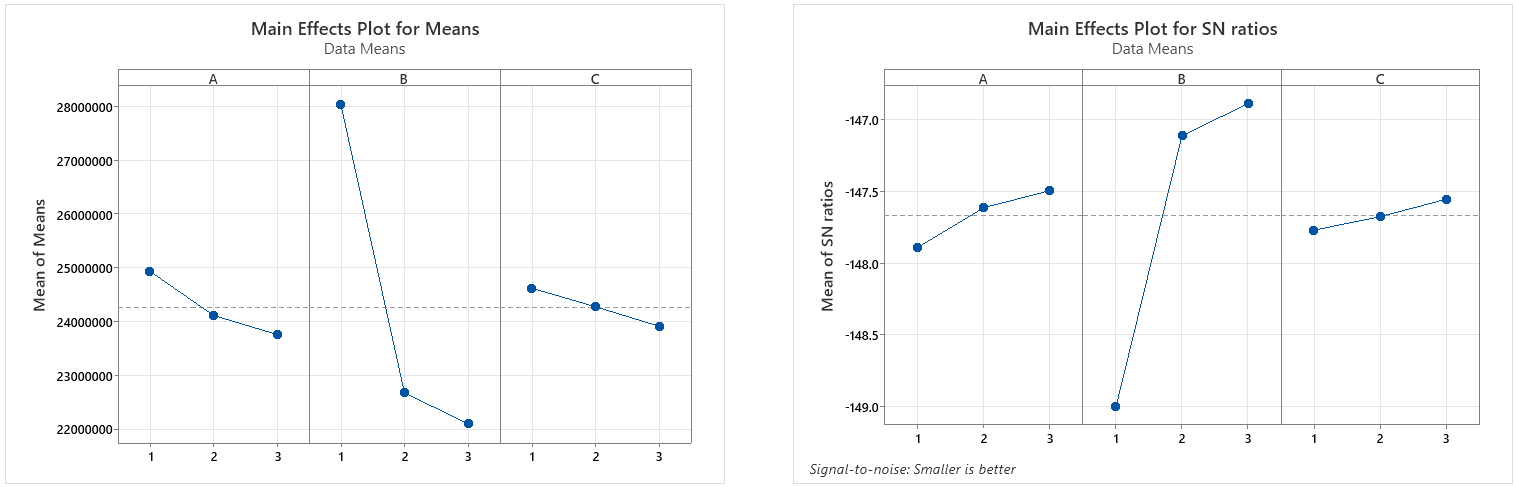}
\caption{The main effects plot for S⁄N of WOA.}
\label{fig:taguchi-woa}
\end{figure}

In addition, Figures~\ref{fig:taguchi-gwo} and~\ref{fig:taguchi-woa} are the main effects plots of S/N ratios.

\subsection{Performance analysis}
A detailed analysis of algorithms is provided in this section. 15 numerical examples with different dimensions (number of products and retailers) are solved by the algorithms. Each algorithm is run in 10 replications for each generated numerical example. We consider four measures, including average RPD, average RDI, average CPU time, and average standard deviation, to asses and compare the efficiency of metaheuristics. The average RPD and RDI as two relative error indicators are used to investigate the solution's quality. These measures are calculated based on the following equations:
\begin{equation} \label{eq:rpd}
    RPD = \frac{Cur_{sol} - Best_{sol}}{Best_{sol}}
\end{equation}
\begin{equation} \label{eq:rdi}
    RDI = \frac{Cur_{sol} - Best_{sol}}{Worst_{sol} - Best_{sol}}
\end{equation}

where $Cur_{sol}$ is the current obtained solution, $Best_{sol}$ and $Worst_{sol}$ are the best and worst calculated solutions, respectively. Note that the lower values of these measures demonstrate better performance. We utilized average CPU time to demonstrate the required time of algorithms to solve the problem. In addition, the standard deviation measure is used to evaluate the robustness of algorithms in different runs. Table~\ref{tab:results-comparison} summarizes the obtained results of the metaheuristic algorithm.

\begin{table}[H]
\caption{The computational results of the metaheuristic algorithms}
\label{tab:results-comparison}
\centering
\begin{adjustbox}{width=\textwidth}
\begin{tabular}{cclcclcclcclcc}
\hline
\multicolumn{2}{c}{Problem size} &  & \multicolumn{2}{c}{Average RPD} &  & \multicolumn{2}{c}{Average RDI} &  & \multicolumn{2}{c}{Average standard deviation} &  & \multicolumn{2}{c}{Average CPU time} \\ \cline{1-2} \cline{4-5} \cline{7-8} \cline{10-11} \cline{13-14} 
Vendor         & Retailer         &  & GWO            & WOA            &  & GWO            & WOA            &  & GWO                   & WOA                    &  & GWO               & WOA              \\ \hline 
1                & 2             &  & 7.119E-07      & 1.930E-09      &  & 4.218E-01      & 4.275E-01      &  & 0.921                 & 0.003                  &  & 14.569            & 16.158           \\
2                & 3             &  & 5.470E-06      & 2.118E-07      &  & 6.211E-01      & 3.869E-01      &  & 14.125                & 1.059                  &  & 16.919            & 18.362           \\
2                & 6             &  & 3.410E-05      & 1.892E-04      &  & 3.529E-01      & 4.506E-01      &  & 319.281               & 1626.265               &  & 18.117            & 22.751           \\
3                & 3             &  & 7.666E-06      & 3.495E-06      &  & 4.170E-01      & 3.970E-01      &  & 50.032                & 25.751                 &  & 18.228            & 19.260           \\
3                & 5             &  & 3.215E-05      & 2.121E-04      &  & 3.977E-01      & 4.572E-01      &  & 440.535               & 1941.583               &  & 18.992            & 20.257           \\
3                & 6             &  & 4.348E-05      & 5.121E-04      &  & 4.886E-01      & 5.082E-01      &  & 620.624               & 6005.441               &  & 19.902            & 21.817           \\
4                & 3             &  & 6.805E-06      & 3.325E-06      &  & 3.537E-01      & 3.307E-01      &  & 67.446                & 35.478                 &  & 19.843            & 19.148           \\
4                & 4             &  & 3.170E-05      & 2.659E-04      &  & 4.893E-01      & 3.385E-01      &  & 356.253               & 4552.478               &  & 21.917            & 23.415           \\
4                & 5             &  & 5.143E-04      & 1.316E-02      &  & 4.667E-01      & 4.230E-01      &  & 7752.310              & 280268.084             &  & 26.476            & 26.376           \\
5                & 5             &  & 7.856E-05      & 1.055E-03      &  & 3.881E-01      & 3.097E-01      &  & 1868.231              & 23823.369              &  & 22.619            & 23.665           \\
5                & 7             &  & 3.475E-04      & 2.442E-03      &  & 5.998E-01      & 3.487E-01      &  & 5200.742              & 68809.121              &  & 24.730            & 24.539           \\
6                & 4             &  & 3.436E-05      & 6.543E-04      &  & 4.863E-01      & 4.669E-01      &  & 643.454               & 11541.881              &  & 21.909            & 22.957           \\
6                & 5             &  & 2.168E-04      & 1.335E-03      &  & 4.572E-01      & 3.394E-01      &  & 4176.421              & 31895.015              &  & 23.000            & 24.088           \\
6                & 6             &  & 1.950E-04      & 2.059E-03      &  & 6.193E-01      & 4.524E-01      &  & 3689.708              & 48374.604              &  & 25.398            & 24.530           \\
7                & 3             &  & 1.819E-05      & 2.274E-04      &  & 3.379E-01      & 3.774E-01      &  & 325.777               & 4136.163               &  & 21.729            & 22.145           \\ \hline
\multicolumn{2}{c}{Average}      &  & 1.045E-04      & 1.475E-03      &  & 4.598E-01      & 4.009E-01      &  & 1701.724              & 32202.420              &  & 20.957            & 21.965           \\ \hline
\end{tabular}
\end{adjustbox}
\end{table}

\begin{figure}[H]
\centering 
\includegraphics[width=17cm]{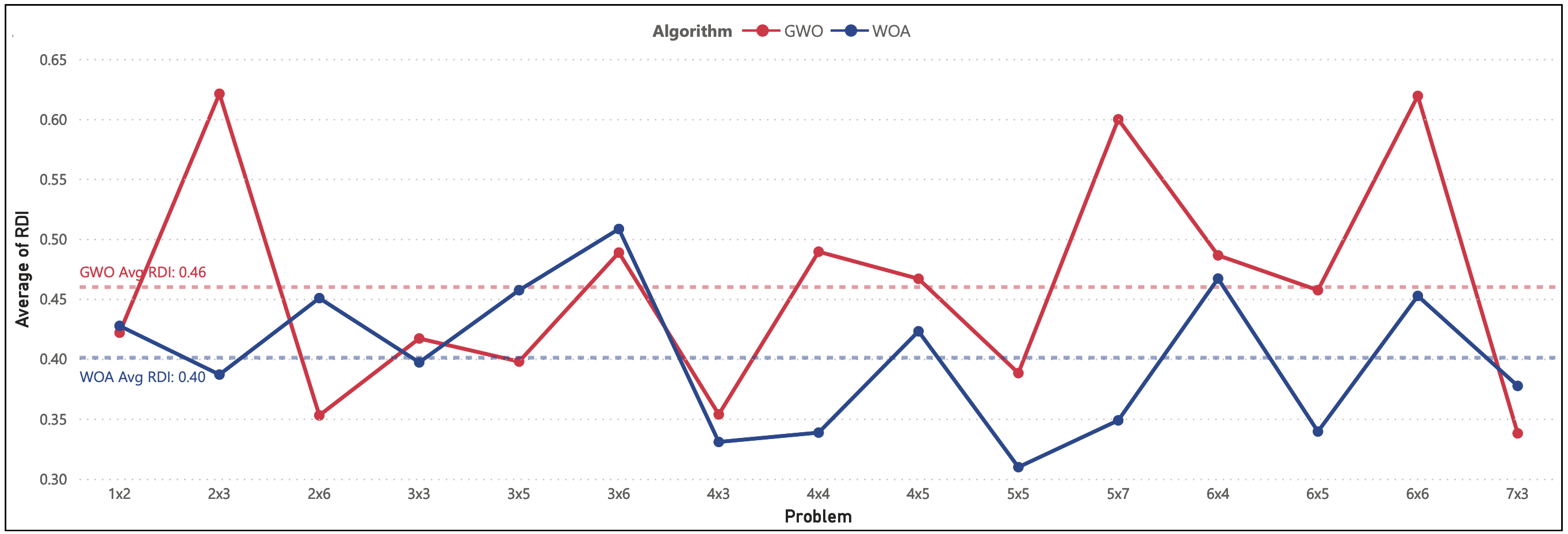}
\caption{The average RDI of metaheuristics for the test examples}
\label{fig:rdi}
\end{figure}

\begin{figure}[H]
\centering 
\includegraphics[width=17cm]{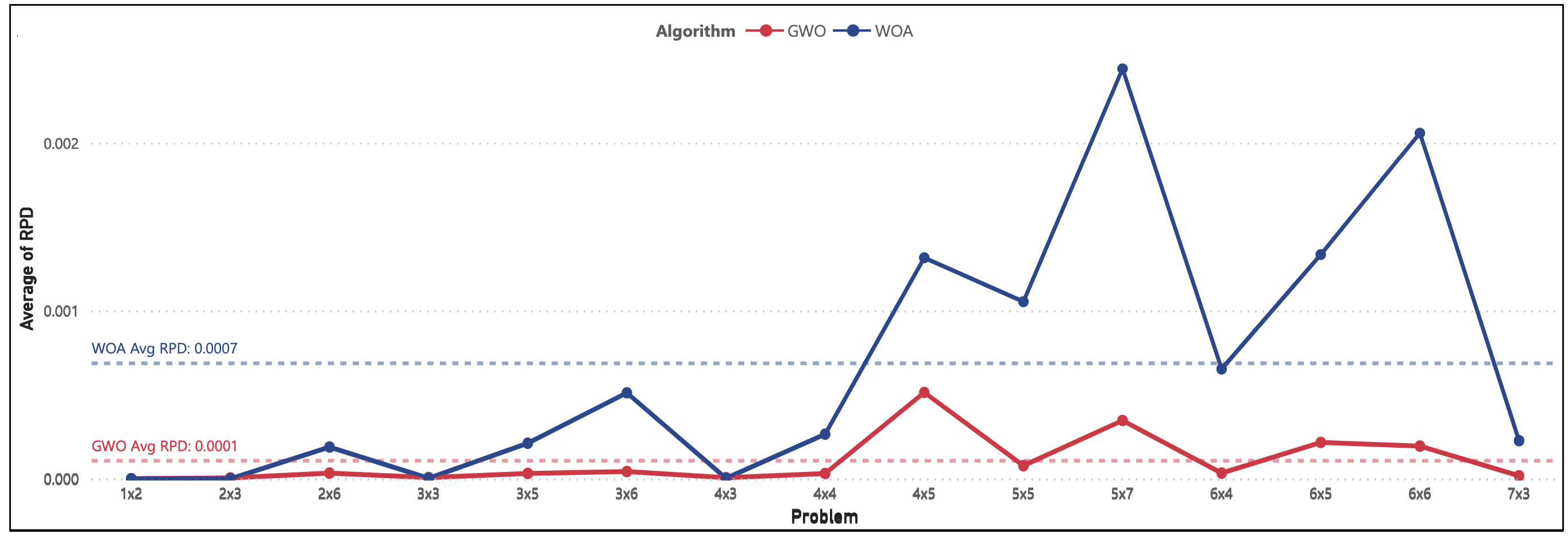}
\caption{The average RPD of metaheuristics for the test examples}
\label{fig:rpd}
\end{figure}

\begin{figure}[H]
\centering 
\includegraphics[width=17cm]{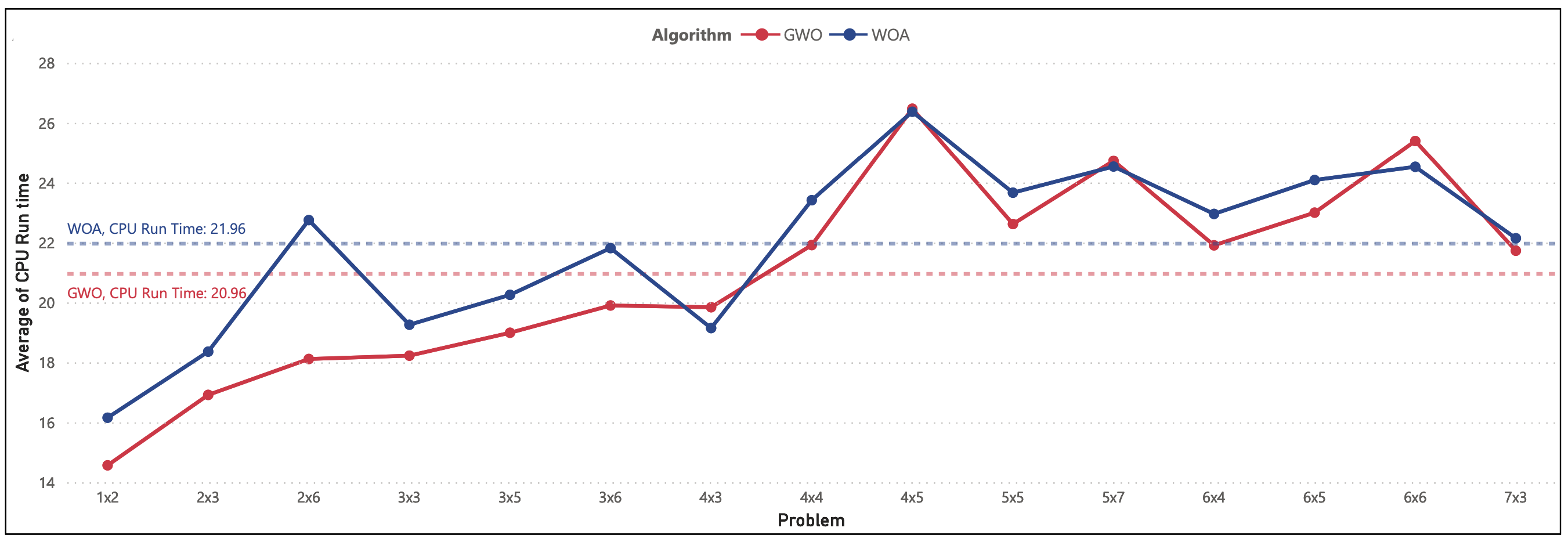}
\caption{The average CPU time of metaheuristics for the test examples.}
\label{fig:cpu}
\end{figure}

\begin{figure}[H]
\centering 
\includegraphics[width=17cm]{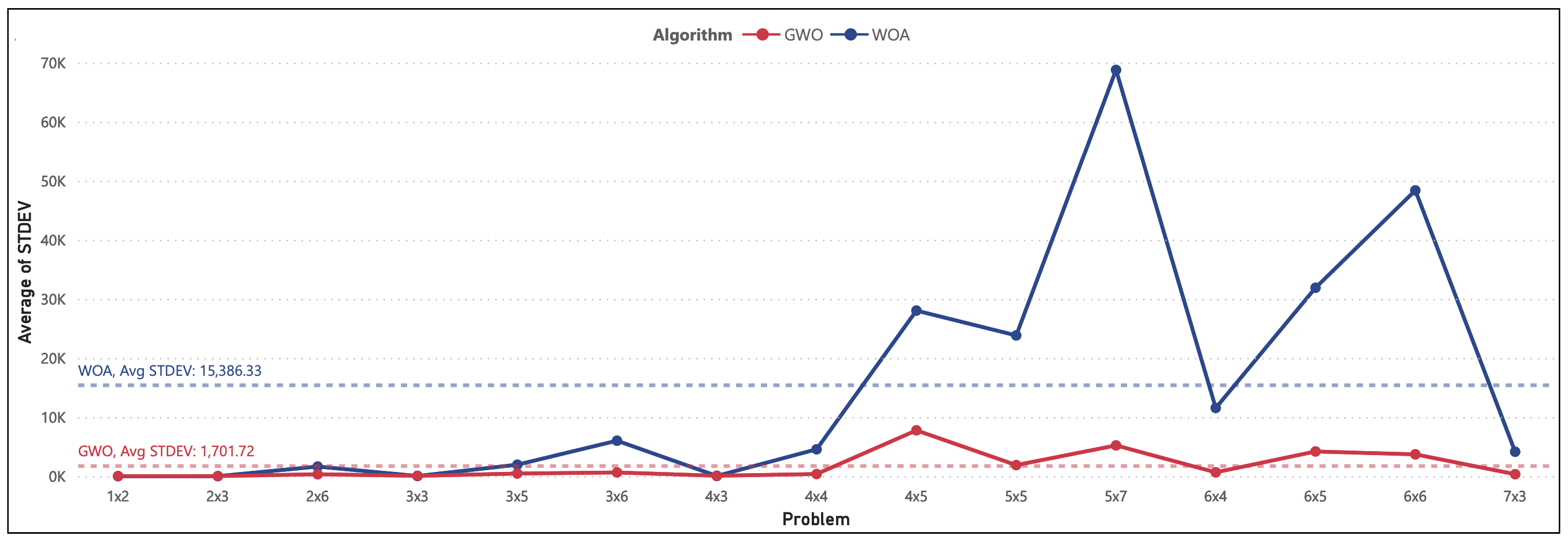}
\caption{The average standard deviation of metaheuristics for the test examples}
\label{fig:stdev}
\end{figure}

Considering the RPD and RDI measures, the algorithms are competitive. As can be seen, GWO reaches lower RPD in most cases. On the other hand, the RDI values of WOA are lower than GWO. The results express that GWO is more robust than WOA. As can be seen, the variation of the calculated solutions by GWO is less than WOA, and GWO has less standard deviation for different numerical examples. Also, the CPU time of the algorithms is very competitive. However, GWO solves most problems in less amount of time. The schematic comparisons of results are presented in Figures \ref{fig:rdi} to \ref{fig:stdev} to provide better insight.

\begin{figure*}[hbt!]
    \centering
    \begin{subfigure}[b]{0.475\textwidth}
        \centering
        \includegraphics[width=\textwidth]{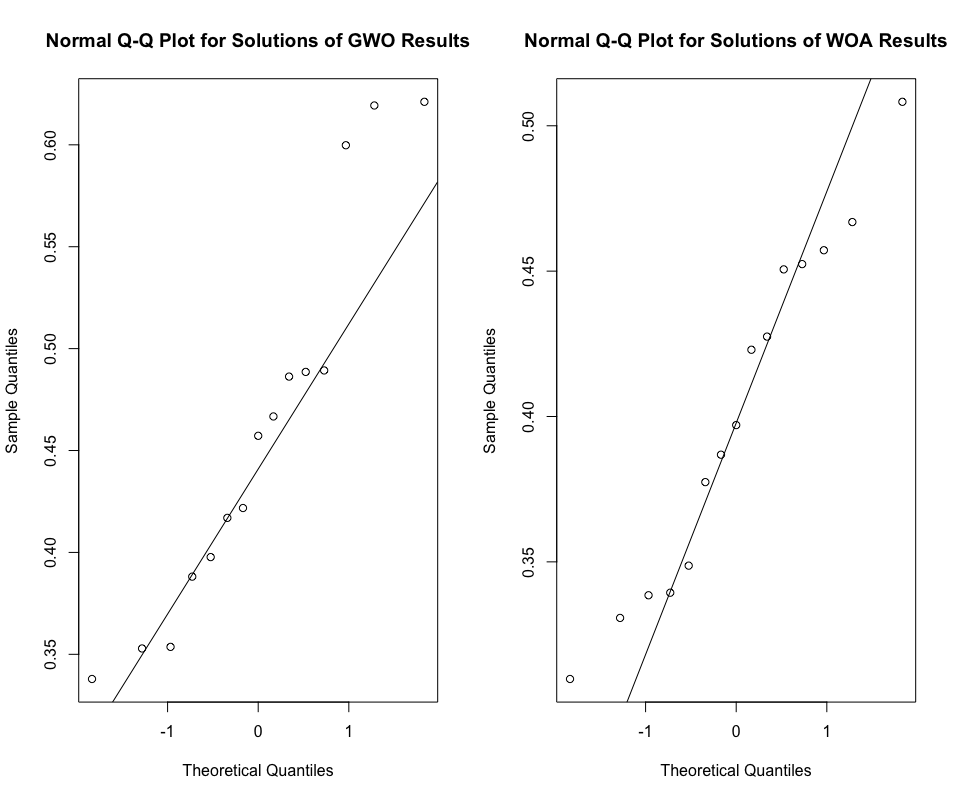}
        \caption[]
        {{\small RDI}}    
        \label{fig:qq-rdi}
    \end{subfigure}
    \hfill
    \begin{subfigure}[b]{0.475\textwidth}  
        \centering 
        \includegraphics[width=\textwidth]{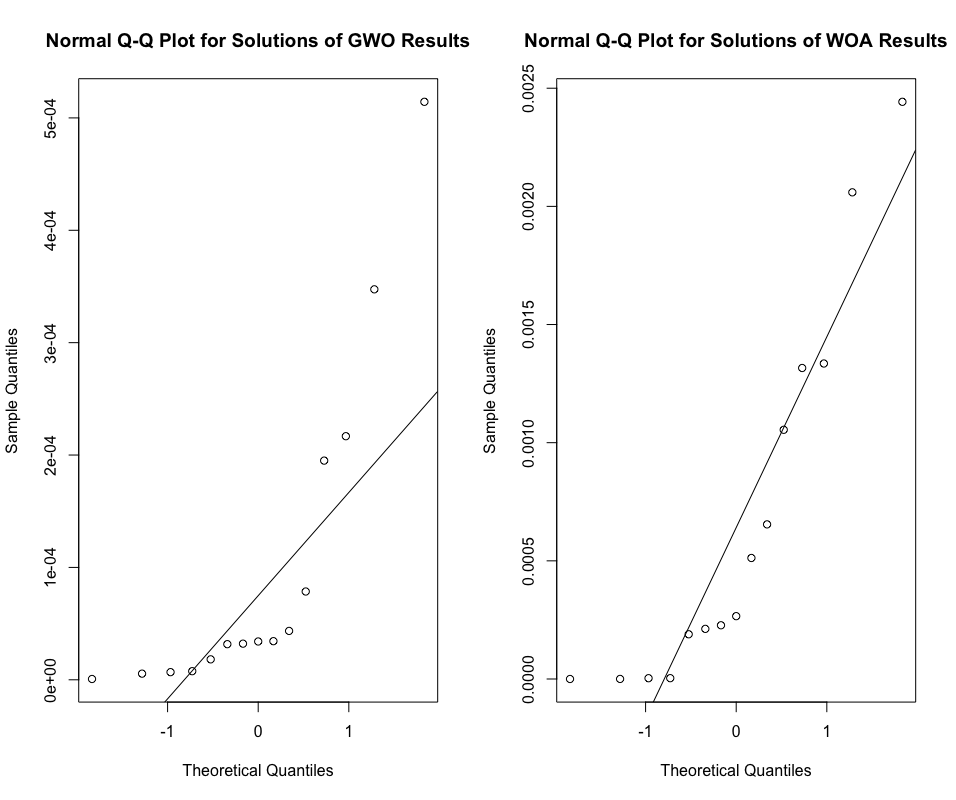}
        \caption[]
        {{\small RPD}}    
        \label{fig:qq-rpd}
    \end{subfigure}
    \vskip\baselineskip
    \begin{subfigure}[b]{0.475\textwidth}   
        \centering 
        \includegraphics[width=\textwidth]{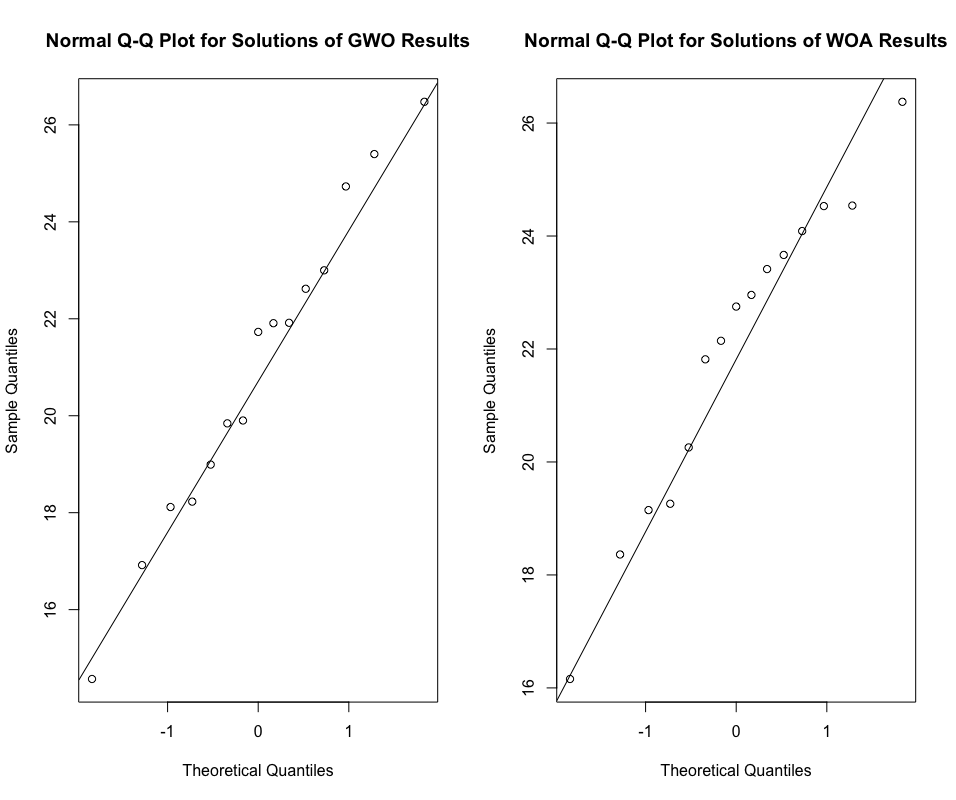}
        \caption[]
        {{\small CPU run time}}    
        \label{fig:qq-cpu}
    \end{subfigure}
    \hfill
    \begin{subfigure}[b]{0.475\textwidth}   
        \centering 
        \includegraphics[width=\textwidth]{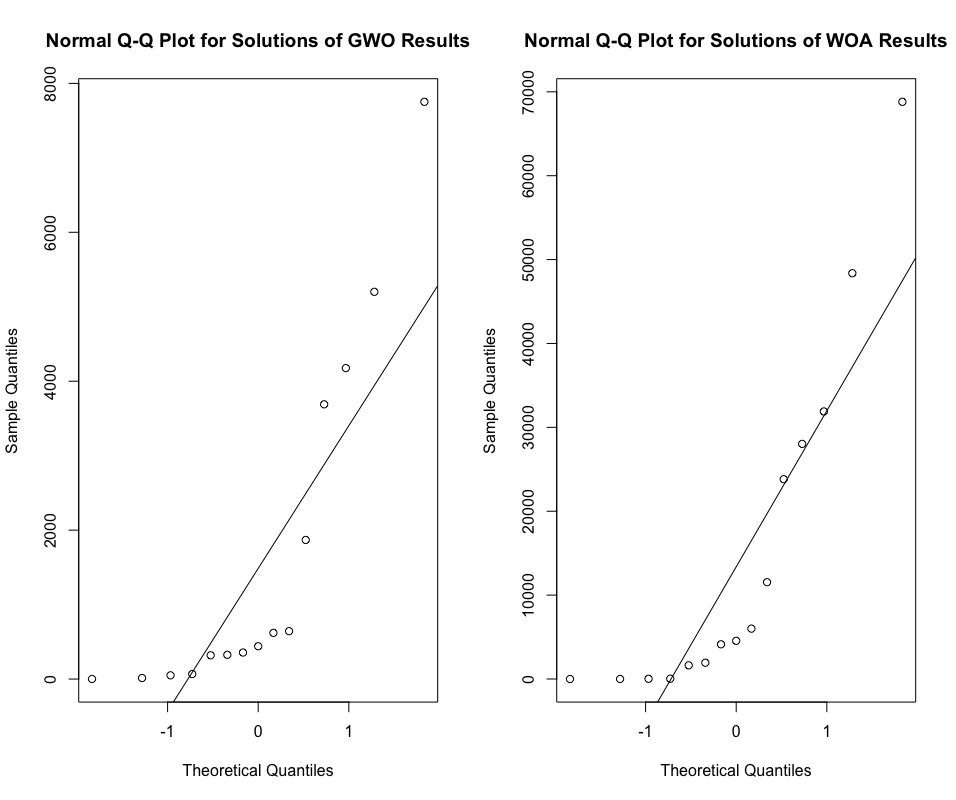}
        \caption[]
        {{\small Standard deviation}}    
        \label{fig:qq-stdev}
    \end{subfigure}
    \caption[]
    {\small The normal Q-Q plots for performance indicators} 
    \label{fig:qq}
\end{figure*}

In continuing, the performance of algorithms is compared statistically. Here, we use statistical hypothesis testing to see whether there is a significant difference between the performance of metaheuristics. All tests and comparisons are performed in $\alpha=0.05$ significance level. The paired-sample t-test and Wilcoxon signed-rank test are parametric and non-parametric tests for the statistical comparison of two populations. To select the proper test, we need to evaluate the normality distribution of considered measures for the solution of algorithms. Therefore, the normal probability plots are provided and presented in Figure ~\ref{fig:qq}. 

Based on the results, we use Wilcoxon signed-rank test for average RPD and average standard deviation and test the difference between average RDI and average CPU time by paired-sample t-test. The results of these tests are provided in Table~\ref{tab:sresults-comparison}.

\begin{table}[H]
\caption{The results of statistical test for algorithms comparison}
\label{tab:sresults-comparison}
\centering
\begin{tabular}{llll}
\hline
Measure & Test & P-value & Superior algorithm \\
\hline
Average RPD & Wilcoxon signed-rank test & 0.002 & GWO \\
Average RDI & Paired sample & 0.048 & WOA \\
Average CPU time & Paired sample & 0.010 & GWO \\
Average Standard deviation & Wilcoxon signed-rank test & 0.002 & GWO \\
\hline
\end{tabular}
\end{table}

\begin{figure*}[hbt!]
    \centering
    \begin{subfigure}[b]{0.475\textwidth}
        \centering
        \includegraphics[width=\textwidth]{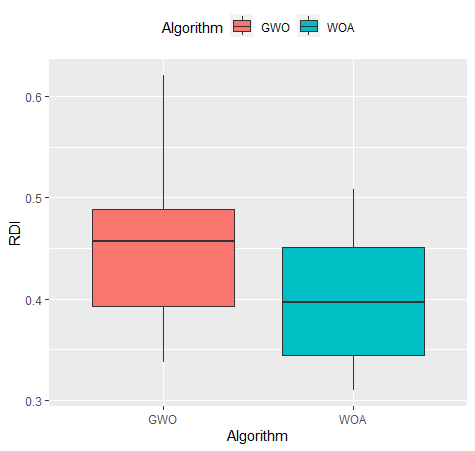}
        \caption[]
        {{\small RDI}}    
        \label{fig:box-rdi}
    \end{subfigure}
    \hfill
    \begin{subfigure}[b]{0.475\textwidth}  
        \centering 
        \includegraphics[width=\textwidth]{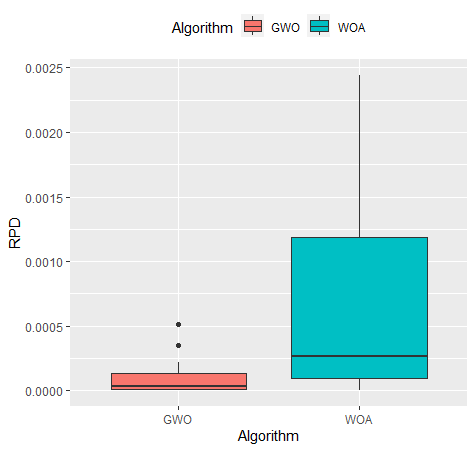}
        \caption[]
        {{\small RPD}}    
        \label{fig:box-rpd}
    \end{subfigure}
    \vskip\baselineskip
    \begin{subfigure}[b]{0.475\textwidth}   
        \centering 
        \includegraphics[width=\textwidth]{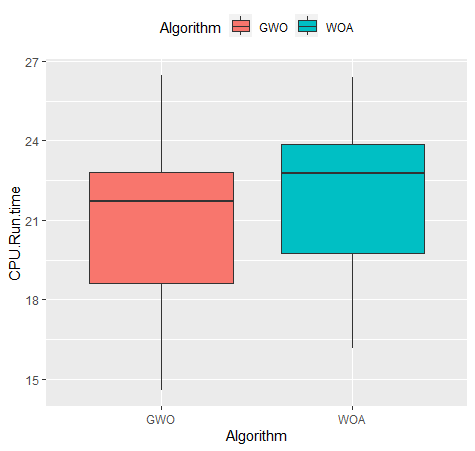}
        \caption[]
        {{\small CPU run time}}    
        \label{fig:box-cpu}
    \end{subfigure}
    \hfill
    \begin{subfigure}[b]{0.475\textwidth}   
        \centering 
        \includegraphics[width=\textwidth]{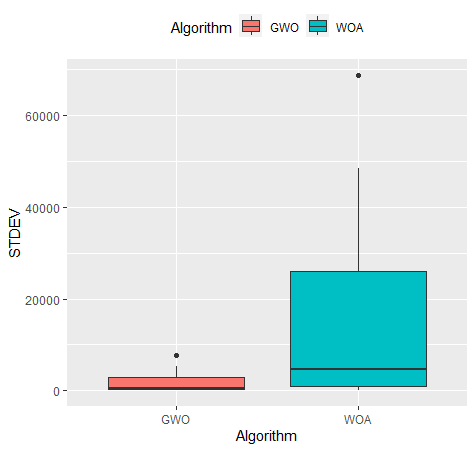}
        \caption[]
        {{\small Standard deviation}}    
        \label{fig:box-stdev}
    \end{subfigure}
    \caption[]
    {\small The boxplots for performance indicators} 
    \label{fig:boxplots}
\end{figure*}

In this table, the p-value of tests is less than 0.05. Consequently, we can infer the significant difference of metaheuristics at this significance level. Considering the results, GWO is the superior algorithm regarding average RPD, CPU time, and standard deviation. However, the WOA algorithm has a significantly better performance in terms of average RDI. In addition, the boxplot of each performance measure is also presented in Figure~\ref{fig:boxplots}.

As can be seen, the boxplots of GWO are lower than the boxplots of WOA for all measures except the average RDI.

\subsection{Sensitivity analysis}
\begin{figure*}[hbt!]
    \centering
    \begin{subfigure}[b]{0.475\textwidth}   
        \centering 
        \includegraphics[width=\textwidth]{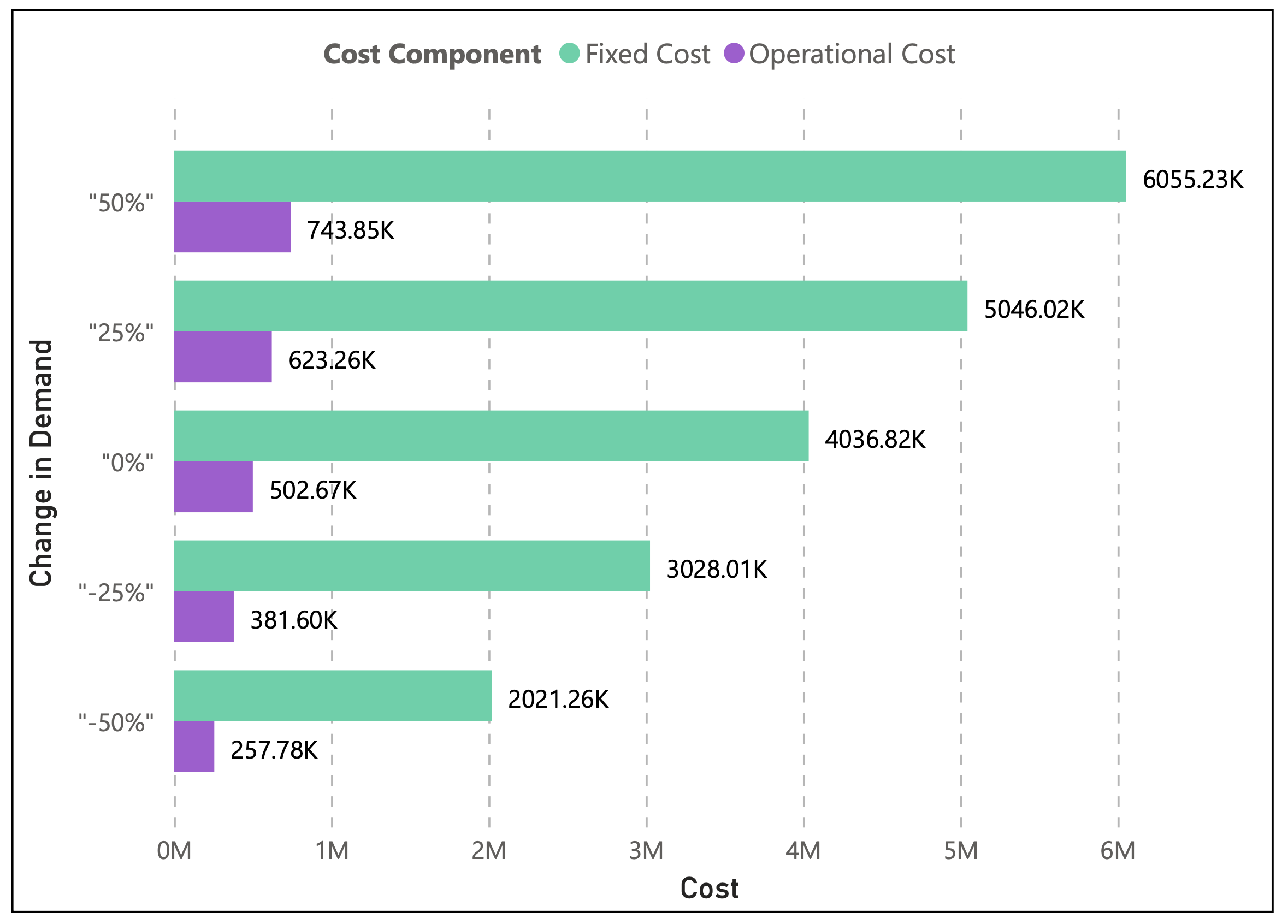}
        \caption[]
        {{\small Fixed and operational cost components}}    
        \label{fig:sen2}
    \end{subfigure}
    \hfill
    \begin{subfigure}[b]{0.475\textwidth}   
        \centering 
        \includegraphics[width=\textwidth]{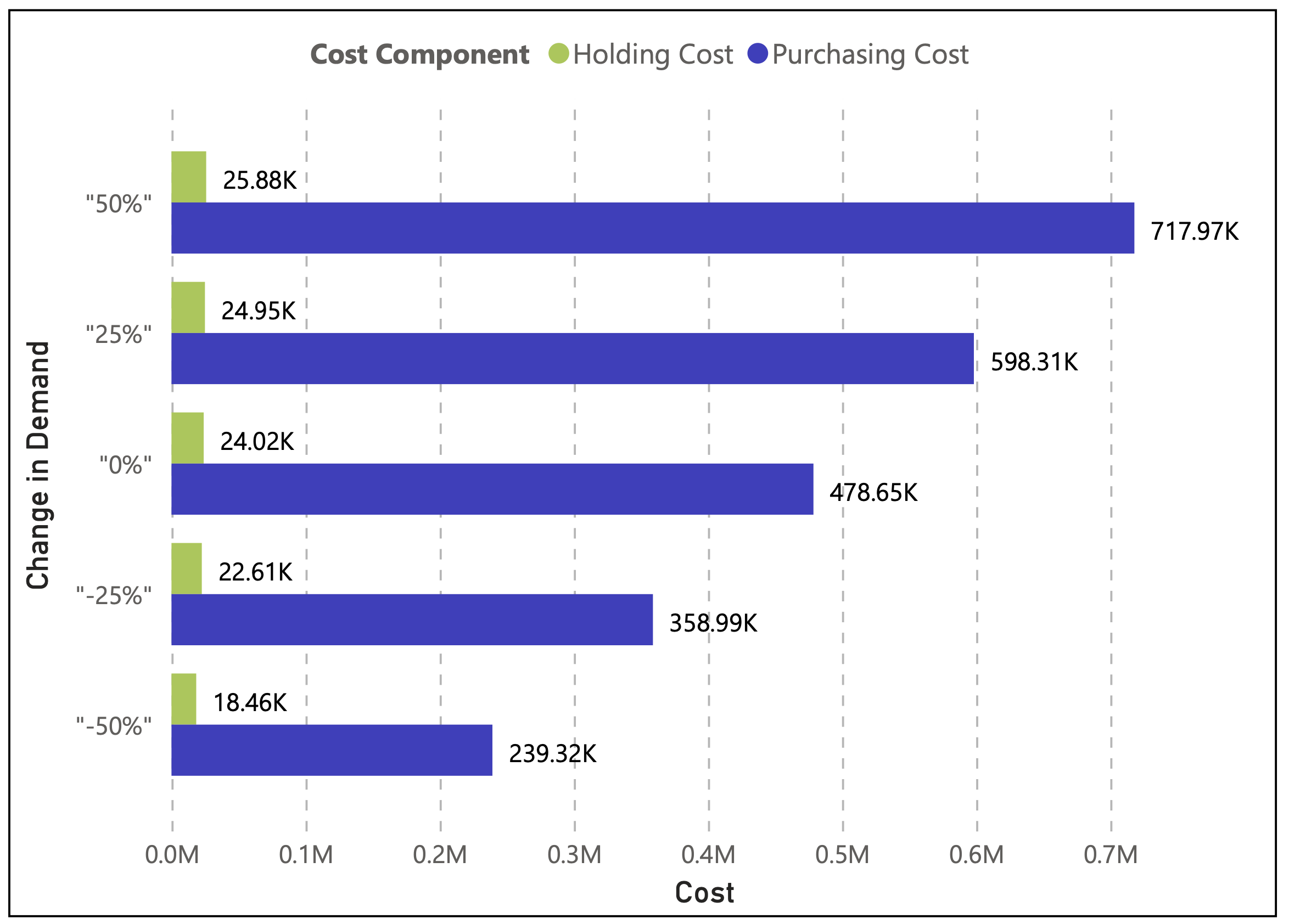}
        \caption[]
        {{\small Operational cost components}}    
        \label{fig:sen3}
    \end{subfigure}
    \caption[]
    {\small Sensitivity Analysis for cost components by changing in demand rate} 
    \label{fig:sensitivity for costs}
\end{figure*}

Sensitivity analysis is a systematic approach that aims to provide more insights for managers into the system's considering the variability of parameters \cite{sadeghi2023mixed}. As the last step, sensitivity analysis is carried out to investigate the impact of change in demand as one of the main parameters on the cost components of the inventory system.  Since the SQP can calculate the optimal results, this is used to perform the sensitivity analysis. For this goal, we consider the change in parameters at -50\% to +50\% rates. In addition, we categorize the cost components to provide a better insight. The fixed costs include the fixed ordering cost of the vendor and retailers and the fixed recovery cost of retailers. Also, the purchasing cost and holding cost of retailers are the operational costs. Figure~\ref{fig:sensitivity for costs} shows the details of obtained results.

As evident, the increasing demand negatively impacts the total cost of the system. In fact, more demand satisfaction requires more ordering, holding, and recovery of products. Based on the results, the fixed cost parameters are more sensitive than the operational cost parameters to changes in demand. In addition, the purchasing cost is more impacted than the holding cost, when the demand is violated. This means managers should focus more on controlling and reducing fixed components to improve the system's performance. Various approaches, such as investment in the infrastructures of the supply chain, can help the managers with this goal.

\section{Conclusion}\label{sec:con}
In this research, we developed a  multi-product model for inventory management of reusable items in  a single-vendor multi-retailer two-level supply chain. The presented problem considered the stochastic limitation on the system's resources, such as the maximum budget, storage capacity, etc. This uncertainty was handled by the chance-constrained programming approach. The problem aimed to establish the optimal ordering and recovery of each product for each retailer so that the overall cost of the supply chain is minimized. Regarding the nonlinearity of the constrained model, GWO and WOA algorithms, as two novel metaheuristics, were proposed as the solution approach. Moreover, the SQP exact algorithm is proposed to assess the performance of GWO and WOA and further analysis. The parameters of metaheuristics are tuned by the Taguchi statistical method.

The computational results of a small numerical example show the power of metaheuristic algorithms in searching the solution space. The obtained results revealed that the difference between the solutions of the SQP exact method and metaheuristic algorithm is not significant. The performance of metaheuristic algorithms was extensively compared using 15 different sizes of numerical examples. The results showed that GWO is preferable to WOA for the presented problem. More specifically, this algorithm obtains higher-quality solutions in less amount of time. We found that GWO is more robust, as the variation of the algorithm solutions is lower than WOA in different runs. Finally, the impact of variation in the demand parameter was evaluated by analyzing sensitivity. The results showed that the fixed cost components are more sensitive than the operational cost components to changes in demand for products.

For future research, uncertainty in other parameters can be considered, and research can model the problem using the stochastic or fuzzy programming method. Furthermore, multi-criteria decision-making methods, such as ABC analysis, can be employed to classify products and improve the performance of the system. Other heuristic and metaheuristic algorithms can be used as solution approaches, and the author can compare them with GWO and WCA.


\newpage
\bibliographystyle{elsarticle-num}
\bibliography{cas-refs}

\begin{thebibliography}{10}
\expandafter\ifx\csname url\endcsname\relax
  \def\url#1{\texttt{#1}}\fi
\expandafter\ifx\csname urlprefix\endcsname\relax\def\urlprefix{URL }\fi
\expandafter\ifx\csname href\endcsname\relax
  \def\href#1#2{#2} \def\path#1{#1}\fi

\bibitem{chang2016supply}
W.~Chang, A.~E. Ellinger, K.~K. Kim, G.~R. Franke, Supply chain integration and
  firm financial performance: A meta-analysis of positional advantage mediation
  and moderating factors, European Management Journal 34~(3) (2016) 282--295.

\bibitem{taghiyeh2020multi}
S.~Taghiyeh, D.~C. Lengacher, A.~H. Sadeghi, A.~Sahebifakhrabad, R.~B.
  Handfield, A multi-phase approach for product hierarchy forecasting in supply
  chain management: application to monarchfx inc, arXiv preprint
  arXiv:2006.08931 (2020).

\bibitem{fakhrabad2023impact}
A.~S. Fakhrabad, A.~H. Sadeghi, E.~Kemahlioglu-Ziya, R.~B. Handfield,
  H.~Tohidi, I.~V. Farahani, The impact of opioid prescribing limits on drug
  usage in south carolina: A novel geospatial and time series data analysis,
  arXiv preprint arXiv:2301.08878 (2023).

\bibitem{nikoubin2023relax}
A.~Nikoubin, M.~Mahnam, G.~Moslehi, A relax-and-fix pareto-based algorithm for
  a bi-objective vaccine distribution network considering a mix-and-match
  strategy in pandemics, Applied Soft Computing 132 (2023) 109862.

\bibitem{sahebi2023evaluating}
A.~Sahebi-Fakhrabad, A.~H. Sadeghi, R.~Handfield, Evaluating state-level
  prescription drug monitoring program (pdmp) and pill mill effects on opioid
  consumption in pharmaceutical supply chain, in: Healthcare, Vol.~11,
  Multidisciplinary Digital Publishing Institute, 2023, p. 437.

\bibitem{hasan2023multi}
K.~W. Hasan, S.~M. Ali, S.~K. Paul, G.~Kabir, Multi-objective closed-loop green
  supply chain model with disruption risk, Applied Soft Computing (2023)
  110074.

\bibitem{braglia2003modelling}
M.~Braglia, L.~Zavanella, Modelling an industrial strategy for inventory
  management in supply chains: The'consignment stock'case, International
  Journal of Production Research 41~(16) (2003) 3793--3808.

\bibitem{disney2003effect}
S.~M. Disney, D.~R. Towill, The effect of vendor managed inventory (vmi)
  dynamics on the bullwhip effect in supply chains, International journal of
  production economics 85~(2) (2003) 199--215.

\bibitem{harris1990many}
F.~W. Harris, How many parts to make at once, Operations research 38~(6) (1990)
  947--950.

\bibitem{mokhtari2018joint}
H.~Mokhtari, Joint ordering and reuse policy for reusable items inventory
  management, Sustainable Production and Consumption 15 (2018) 163--172.

\bibitem{fallahi2022constrained}
A.~Fallahi, E.~A. Bani, S.~T.~A. Niaki, A constrained multi-item eoq inventory
  model for reusable items: Reinforcement learning-based differential evolution
  and particle swarm optimization, Expert Systems with Applications 207 (2022)
  118018.

\bibitem{pasandideh2015optimization}
S.~H.~R. Pasandideh, S.~T.~A. Niaki, A.~Gharaei, Optimization of a multiproduct
  economic production quantity problem with stochastic constraints using
  sequential quadratic programming, Knowledge-Based Systems 84 (2015) 98--107.

\bibitem{khalilpourazari2019robust}
S.~Khalilpourazari, S.~H.~R. Pasandideh, A.~Ghodratnama, Robust possibilistic
  programming for multi-item eoq model with defective supply batches: Whale
  optimization and water cycle algorithms, Neural computing and applications 31
  (2019) 6587--6614.

\bibitem{taft1918most}
E.~Taft, The most economical production lot, Iron Age 101~(18) (1918)
  1410--1412.

\bibitem{covert1973eoq}
R.~P. Covert, G.~C. Philip, An eoq model for items with weibull distribution
  deterioration, AIIE transactions 5~(4) (1973) 323--326.

\bibitem{mokhtari2022economic}
H.~Mokhtari, A.~Salmasnia, A.~Fallahi, Economic production quantity under
  possible substitution: A scenario analysis approach, International Journal of
  Industrial Engineering \& Production Research 33~(1) (2022) 1--17.

\bibitem{drezner1995eoq}
Z.~Drezner, H.~Gurnani, B.~A. Pasternack, An eoq model with substitutions
  between products, Journal of the Operational Research Society 46~(7) (1995)
  887--891.

\bibitem{rezaei2014economic}
J.~Rezaei, Economic order quantity for growing items, International Journal of
  Production Economics 155 (2014) 109--113.

\bibitem{mokhtari2020extended}
H.~Mokhtari, J.~Asadkhani, Extended economic production quantity models with
  preventive maintenance, Scientia Iranica 27~(6) (2020) 3253--3264.

\bibitem{porteus1985investing}
E.~L. Porteus, Investing in reduced setups in the eoq model, Management Science
  31~(8) (1985) 998--1010.

\bibitem{taleizadeh2013eoq}
A.~A. Taleizadeh, D.~W. Pentico, M.~S. Jabalameli, M.~Aryanezhad, An eoq model
  with partial delayed payment and partial backordering, Omega 41~(2) (2013)
  354--368.

\bibitem{chung2003optimal}
K.-J. Chung, Y.-F. Huang, The optimal cycle time for epq inventory model under
  permissible delay in payments, International Journal of Production Economics
  84~(3) (2003) 307--318.

\bibitem{chen2013carbon}
X.~Chen, S.~Benjaafar, A.~Elomri, The carbon-constrained eoq, Operations
  Research Letters 41~(2) (2013) 172--179.

\bibitem{taleizadeh2018sustainable}
A.~A. Taleizadeh, V.~R. Soleymanfar, K.~Govindan, Sustainable economic
  production quantity models for inventory systems with shortage, Journal of
  cleaner production 174 (2018) 1011--1020.

\bibitem{tripathy2003eoq}
P.~Tripathy, W.~Wee, P.~R. Majhi, An eoq model with process reliability
  considerations, Journal of the Operational Research Society 54~(5) (2003)
  549--554.

\bibitem{cheng1990eoq}
T.~Cheng, An eoq model with pricing consideration, Computers \& industrial
  engineering 18~(4) (1990) 529--534.

\bibitem{subramanyam1981eoq}
E.~S. Subramanyam, S.~Kumaraswamy, Eoq formula under varying marketing policies
  and conditions, AIIE Transactions 13~(4) (1981) 312--314.

\bibitem{salameh2000economic}
M.~K. Salameh, M.~Y. Jaber, Economic production quantity model for items with
  imperfect quality, International journal of production economics 64~(1-3)
  (2000) 59--64.

\bibitem{chandra1985effects}
M.~J. Chandra, M.~L. Bahner, The effects of inflation and the time value of
  money on some inventory systems, International Journal of Production Research
  23~(4) (1985) 723--730.

\bibitem{fallahi2021sustainable}
A.~Fallahi, M.~Azimi-Dastgerdi, H.~Mokhtari, A sustainable production-inventory
  model joint with preventive maintenance and multiple shipments for imperfect
  quality items, Scientia Iranica (2021).

\bibitem{khan2011economic}
M.~Khan, M.~Y. Jaber, M.~Bonney, An economic order quantity (eoq) for items
  with imperfect quality and inspection errors, International Journal of
  Production Economics 133~(1) (2011) 113--118.

\bibitem{pasandideh2011genetic}
S.~H.~R. Pasandideh, S.~T.~A. Niaki, A.~R. Nia, A genetic algorithm for vendor
  managed inventory control system of multi-product multi-constraint economic
  order quantity model, Expert Systems with Applications 38~(3) (2011)
  2708--2716.

\bibitem{pasandideh2014optimization}
S.~H.~R. Pasandideh, S.~T.~A. Niaki, M.~Hemmati~Far, Optimization of vendor
  managed inventory of multiproduct epq model with multiple constraints using
  genetic algorithm, The International Journal of Advanced Manufacturing
  Technology 71 (2014) 365--376.

\bibitem{pasandideh2011parameter}
S.~H.~R. Pasandideh, S.~T.~A. Niaki, N.~Tokhmehchi, A parameter-tuned genetic
  algorithm to optimize two-echelon continuous review inventory systems, Expert
  Systems with Applications 38~(9) (2011) 11708--11714.

\bibitem{taleizadeh2011multiple}
A.~A. Taleizadeh, S.~T.~A. Niaki, F.~Barzinpour, Multiple-buyer multiple-vendor
  multi-product multi-constraint supply chain problem with stochastic demand
  and variable lead-time: a harmony search algorithm, Applied Mathematics and
  Computation 217~(22) (2011) 9234--9253.

\bibitem{chen2014retailer}
S.-C. Chen, L.~E. C{\'a}rdenas-Barr{\'o}n, J.-T. Teng, Retailer’s economic
  order quantity when the supplier offers conditionally permissible delay in
  payments link to order quantity, International journal of production
  Economics 155 (2014) 284--291.

\bibitem{cardenas2015multi}
L.~E. C{\'a}rdenas-Barr{\'o}n, S.~S. Sana, Multi-item eoq inventory model in a
  two-layer supply chain while demand varies with promotional effort, Applied
  Mathematical Modelling 39~(21) (2015) 6725--6737.

\bibitem{khan2017learning}
M.~Khan, M.~Hussain, L.~E. C{\'a}rdenas-Barr{\'o}n, Learning and screening
  errors in an epq inventory model for supply chains with stochastic lead time
  demands, International Journal of Production Research 55~(16) (2017)
  4816--4832.

\bibitem{tiwari2018joint}
S.~Tiwari, L.~E. C{\'a}rdenas-Barr{\'o}n, M.~Goh, A.~A. Shaikh, Joint pricing
  and inventory model for deteriorating items with expiration dates and partial
  backlogging under two-level partial trade credits in supply chain,
  International Journal of Production Economics 200 (2018) 16--36.

\bibitem{karimian2020geometric}
Y.~Karimian, A.~Mirzazadeh, S.~H. Pasandideh, M.~Namakshenas, A geometric
  programming approach for a vendor managed inventory of a multiretailer
  multi-item epq model, RAIRO-Operations Research 54~(5) (2020) 1401--1418.

\bibitem{pourmohammad2021food}
N.~Pourmohammad-Zia, B.~Karimi, J.~Rezaei, Food supply chain coordination for
  growing items: A trade-off between market coverage and cost-efficiency,
  International Journal of Production Economics 242 (2021) 108289.

\bibitem{pourmohammad2021dynamic}
N.~Pourmohammad-Zia, B.~Karimi, J.~Rezaei, Dynamic pricing and inventory
  control policies in a food supply chain of growing and deteriorating items,
  Annals of Operations Research (2021) 1--40.

\bibitem{asadkhani2022sustainable}
J.~Asadkhani, A.~Fallahi, H.~Mokhtari, A sustainable supply chain under vmi-cs
  agreement with withdrawal policies for imperfect items, Journal of Cleaner
  Production 376 (2022) 134098.

\bibitem{hill1997single}
R.~M. Hill, The single-vendor single-buyer integrated production-inventory
  model with a generalised policy, European journal of operational research
  97~(3) (1997) 493--499.

\bibitem{moheb2023reverse}
H.~Moheb-Alizadeh, A.~H. Sadeghi, M.~K. Jaunich, E.~Kemahlioglu-Ziya, R.~B.
  Handfield, et~al., Reverse logistics network design to estimate the economic
  and environmental impacts of take-back legislation: A case study for e-waste
  management system in washington state, arXiv preprint arXiv:2301.09792
  (2023).

\bibitem{khalilpourazari2016optimization}
S.~Khalilpourazari, S.~H.~R. Pasandideh, S.~T.~A. Niaki, Optimization of
  multi-product economic production quantity model with partial backordering
  and physical constraints: Sqp, sfs, sa, and wca, Applied Soft Computing 49
  (2016) 770--791.

\bibitem{mirjalili2014grey}
S.~Mirjalili, S.~M. Mirjalili, A.~Lewis, Grey wolf optimizer, Advances in
  engineering software 69 (2014) 46--61.

\bibitem{faris2018grey}
H.~Faris, I.~Aljarah, M.~A. Al-Betar, S.~Mirjalili, Grey wolf optimizer: a
  review of recent variants and applications, Neural computing and applications
  30 (2018) 413--435.

\bibitem{jiang2022dsgwo}
J.~Jiang, Z.~Zhao, Y.~Liu, W.~Li, H.~Wang, Dsgwo: An improved grey wolf
  optimizer with diversity enhanced strategy based on group-stage competition
  and balance mechanisms, Knowledge-Based Systems 250 (2022) 109100.

\bibitem{mirjalili2016whale}
S.~Mirjalili, A.~Lewis, The whale optimization algorithm, Advances in
  engineering software 95 (2016) 51--67.

\bibitem{zhao2022multipopulation}
F.~Zhao, H.~Bao, L.~Wang, J.~Cao, J.~Tang, et~al., A multipopulation
  cooperative coevolutionary whale optimization algorithm with a two-stage
  orthogonal learning mechanism, Knowledge-Based Systems 246 (2022) 108664.

\bibitem{chakraborty2021enhanced}
S.~Chakraborty, A.~K. Saha, R.~Chakraborty, M.~Saha, An enhanced whale
  optimization algorithm for large scale optimization problems, Knowledge-Based
  Systems 233 (2021) 107543.

\bibitem{gharehchopogh2019comprehensive}
F.~S. Gharehchopogh, H.~Gholizadeh, A comprehensive survey: Whale optimization
  algorithm and its applications, Swarm and Evolutionary Computation 48 (2019)
  1--24.

\bibitem{gill2011sequential}
P.~E. Gill, E.~Wong, Sequential quadratic programming methods, in: Mixed
  integer nonlinear programming, Springer, 2011, pp. 147--224.

\bibitem{gurwitz1989sequential}
C.~B. Gurwitz, M.~L. Overton, Sequential quadratic programming methods based on
  approximating a projected hessian matrix, SIAM Journal on Scientific and
  Statistical Computing 10~(4) (1989) 631--653.

\bibitem{rockafellar1970convex}
R.~T. Rockafellar, Convex analysis, Vol.~18, Princeton university press, 1970.

\bibitem{bartlett2000active}
R.~Bartlett, A.~Wachter, L.~Biegler, Active set vs. interior point strategies
  for model predictive control, in: Proceedings of the 2000 American Control
  Conference. ACC (IEEE Cat. No. 00CH36334), Vol.~6, IEEE, 2000, pp.
  4229--4233.

\bibitem{khalilpourazari2019modeling}
S.~Khalilpourazari, S.~H.~R. Pasandideh, Modeling and optimization of
  multi-item multi-constrained eoq model for growing items, Knowledge-Based
  Systems 164 (2019) 150--162.

\bibitem{mousavi2021intelligent}
S.~M. Mousavi, S.~Abdullah, S.~T.~A. Niaki, S.~Banihashemi, An intelligent
  hybrid classification algorithm integrating fuzzy rule-based extraction and
  harmony search optimization: Medical diagnosis applications, Knowledge-Based
  Systems 220 (2021) 106943.

\bibitem{roy2010primer}
R.~K. Roy, A primer on the Taguchi method, Society of Manufacturing Engineers,
  2010.

\bibitem{sadeghi2023mixed}
A.~H. Sadeghi, Z.~Sun, A.~S. Fakhrabad, H.~Arzani, R.~B. Handfield, A
  mixed-integer linear formulation for dynamic modified stochastic p-median
  problem in a competitive supply chain network design, arXiv preprint
  arXiv:2301.11502 (2023).

\end{thebibliography}

\end{document}